\documentclass[12pt]{article}
\usepackage{amsmath,amssymb,amsfonts,amscd,amsxtra}
\usepackage{latexsym}

\input{epsf.sty}
\usepackage{caption}
\usepackage{epsfig}
\usepackage{graphicx,color}
\usepackage{graphics}
\usepackage{subfigure}
\usepackage{appendix}
\usepackage{algorithm}
\usepackage{array}
\usepackage{booktabs}
\usepackage{algpseudocode}
\usepackage[font=footnotesize]{caption}
\usepackage[top=1in, bottom=1in, left=1in, right=1in]{geometry}

\renewcommand{\theequation}{\thesection.\arabic{equation}}

\newtheorem{thm}{Theorem}[section]

\newtheorem{lem}[thm]{Lemma}
\newtheorem{prop}[thm]{Proposition}

\renewcommand{\Im}{{\mbox{Im}}}
\renewcommand{\Re}{{\mbox{Re}}}

\newcommand{\norm}[1]{\left\Vert#1\right\Vert}

\newcommand{\qed}{\hfill \ensuremath{\square}}

\renewcommand\appendix{\par
  \setcounter{section}{0}
  \setcounter{subsection}{0}
  \setcounter{figure}{0}
  \setcounter{table}{0}
  \renewcommand\thesection{Appendix \Alph{section}}
  \renewcommand\theequation{\Alph{section}.\arabic{equation}}
  \renewcommand\thefigure{\Alph{section}.\arabic{figure}}
  \renewcommand\thetable{\Alph{section}.\arabic{table}}
  \renewcommand\thethm{\Alph{section}.\arabic{thm}}
}

\numberwithin{equation}{section}

\date{}

\title{An adaptive boundary element method for the transmission problem with hyperbolic metamaterials}
\author{
Junshan Lin \thanks{\footnotesize Department of Mathematics and Statistics, Auburn University, Auburn, AL 36849 (jzl0097@
auburn.edu). Junshan Lin was partially supported by the NSF grant DMS-2011148.} }

\begin{document}
\maketitle

%% ------------------------------------------------------------------
%% ABSTRACT
%% ------------------------------------------------------------------
\begin{abstract}
In this work we present an adaptive boundary element method for computing the electromagnetic response of wave interactions in hyperbolic metamaterials.
One unique feature of hyperbolic metamaterial is the strongly directional wave in its propagating cone, which induces sharp transition for the solution of the integral equation 
across the cone boundary when wave starts to decay or grow exponentially. In order to avoid a global refined mesh over the whole boundary,
we employ a two-level a posteriori error estimator and an adaptive mesh refinement procedure to resolve the singularity locally for
the solution of the integral equation. Such an adaptive procedure allows for the reduction of the degree of freedom significantly for the integral equation solver
while achieving desired accuracy for the solution. In addition, to resolve the fast transition of the fundamental solution and its derivatives accurately across the propagation cone boundary,
adaptive numerical quadrature rules are applied to evaluate the integrals for the stiff matrices.
Finally, in order to formulate the integral equations over the boundary, 
we also derive the limits of layer potentials and their derivatives in the hyperbolic media when the target points approach the boundary .
\end{abstract}


% \begin{AMS}
%  35C20,  35Q60, 35P30.
% \end{AMS}

\setcounter{equation}{0}
\setlength{\arraycolsep}{0.25em}
\section{Introduction}\label{sec:introduction}
\subsection{Background}

% \begin{figure}
% \begin{center}
% \includegraphics[height=6cm,width=14cm]{periodic_slit.pdf}
% \vspace*{-20pt}
% \caption{
% }\label{fig:prob_geo}
% \end{center}
% \end{figure}

Hyperbolic metamaterials are a class of anisotropic electromagnetic materials for which one of
the principal components of their relative permittivity or permeability tensors attain opposite sign of the other two principal components:
\begin{equation*}
\boldsymbol{\varepsilon}(x)=\left[
\begin{array}{ccc}
\varepsilon_\parallel(x) & 0 & 0 \\
0  & \varepsilon_\perp(x)  & 0 \\
0 & 0 & \varepsilon_\perp(x)  
\end{array} \right]
\quad \mbox{or} \quad
\boldsymbol{\mu}(x)=\left[
\begin{array}{ccc}
\mu_\parallel(x) & 0 & 0 \\
0  & \mu_\perp(x)  & 0 \\
0 & 0 &  \mu_\perp(x)
\end{array} \right].
\end{equation*}
In the above, the subscripts $\parallel$ and $\perp$ denotes the component parallel and perpendicular to the anisotropy axis respectively,
and there holds 
$$\Re \varepsilon_\perp \cdot  \Re\varepsilon_\parallel <0 \quad \mbox{or}  \quad \Re\mu_\perp \cdot \Re \mu_\parallel <0. $$
Hyperbolic metamaterials can be realized, for instance, by metal–dielectric structures or by embedding arrays of metallic wires in a dielectric matrix 
by restricting free-electron motion to certain directions \cite{Ferrari, Poddubny, Shekhar}. More recently, hexagonal boron nitride (hBN),
$\alpha$-phase molybdenum trioxide ($\alpha$-Mo$O_3$), $\alpha$-phase vanadium pentoxide ($\alpha$-$V_2O_5$) and a few others
emerge as natural hyperbolic materials that attain opposite signs for the in-plane and out-of-plane components of the dielectric tensor \cite{Caldwell, Dai1, Dai2, Novoselov, Taboada}.

Assume that the media is nonmagentic such that the permeability reduces to the unit tensor,
% When the permittivity values $\varepsilon_\parallel$ and $\varepsilon_\perp$ are constant, 
then the dispersion relation for the time-harmonic (with $e^{-i\omega t}$ dependence) Maxwell's equations
\begin{equation*}
\nabla \times E = i \omega \mu_0 \mu H,  \quad \nabla \times H = -i \omega \varepsilon_0\varepsilon E
\end{equation*}
is given by
\begin{equation}\label{eq:dispersion}
\left(k_1^2+k_2^2 + k_3^2 - \varepsilon_\perp k_0^2 \right) \cdot \left(\frac{k_1^2}{\varepsilon_\perp}+\frac{k_2^2+k_3^2}{\varepsilon_\parallel} - k_0^2 \right) = 0,
\end{equation}
where $\varepsilon_0$ and $\mu_0$ are the free-space permittivity or permeability,  $k_0$ is the free-space wavenumber, 
and $k_1$, $k_2$ and $k_3$ are the $x_1$, $x_2$ and $x_3$ components of the wave vector respectively in the Cartesian coordinate.
The first term in \eqref{eq:dispersion} corresponds  a spherical isofrequency surface for the transverse electric (TE) polarized waves, while the second term 
gives rise to a hyperboloidal isofrequency surface for the transverse magnetic (TM) polarized waves
when $\Re \varepsilon_\perp \cdot  \Re\varepsilon_\parallel <0$.
It is seen that the TM waves remain propagating with arbitrarily large wave vectors in the hyperbolic medium, as opposed to evanescent in an isotropic medium.
This unique property leads to many interesting applications of hyperbolic metamaterials ranging from sub-wavelength light manipulation and 
imaging to spontaneous and thermal emission modification \cite{Cortes, Guo, Jacob, Salandrino, Sreekanth}.

\subsection{Problem formulation}
In this paper, we investigate the computation of the electromagnetic response from the wave interactions in the hyperbolic metamaterials.
We focus on the two-dimensional problem when the medium is invariant along the $x_3$ direction and the wave is TM-polarized with the magnetic field $H=(0, 0, u)^T$.
The TE-polarized case is less interesting as it leads to an isotropic problem with the dispersion relation given by the first term of \eqref{eq:dispersion},
and various existing computational methods can be applied to solve the problem.
The Maxwell's equations in the hyperbolic medium for the TM-polarized polarization reduce to the scalar wave equation
\begin{equation}\label{eq:hyperbolic}
\nabla \cdot (A \nabla u)  + k_0^2 u = 0, 
\end{equation}
where the coefficient matrix
\begin{equation}\label{eq:A}
A(x) =  \left[\begin{array}{cc}
\medskip
\frac{1}{\varepsilon_\perp(x)} & 0  \\
0  & \frac{1}{\varepsilon_\parallel(x)}
\end{array} \right].
\end{equation}
In a homogeneous medium, we set $\varepsilon_\perp(x) \equiv \varepsilon_1$ and $\varepsilon_\parallel(x)\equiv\varepsilon_2$,
where the complex-valued permittivities $\varepsilon_1$ and $\varepsilon_2$ satisfy
$$\Re \varepsilon_1 \cdot  \Re\varepsilon_2 <0 \quad \mbox{and} \quad \Im \varepsilon_j>0 \quad (j=1,2), $$
i.e., the hyperbolic medium is lossy.

Assume that a hyperbolic metamaterial with permittivity values $\varepsilon_\perp(x) \equiv \varepsilon_1^{(1)}$ and $\varepsilon_\parallel(x)\equiv\varepsilon_2^{(1)}$
occupies a bounded simply connected domain $\Omega_1$. 
It is placed in an isotropic medium (e.g, vacuum, silicon) or is embedded in another hyperbolic metamaterial with the permittivity values 
$\varepsilon_\perp(x) \equiv \varepsilon_1^{(2)}$ and $\varepsilon_\parallel(x)\equiv\varepsilon_2^{(2)}$, which occupies the region $\Omega_2= \mathbb{R}^2\backslash\bar\Omega$.
When a near-field source is excited in the interior or exterior domain, the magnetic field $u$ satisfies
\begin{equation}\label{eq:trans_prob}
\nabla \cdot (A_j \nabla u)   + k_0^2 u_j = f_j  \quad \mbox{in} \; \Omega_j, \quad j=1,2,
\end{equation}
in which
\begin{equation*}
A_j  =  \left[\begin{array}{cc}
\left(\varepsilon_1^{(j)}\right)^{-1}  & 0  \\
0  & \left(\varepsilon_2^{(j)}\right)^{-1}
\end{array} \right].
\end{equation*}
There holds $\varepsilon_1^{(2)}=\varepsilon_2^{(2)}$ when the exterior region is isotropic. The source $f_j$ attains a compact support in $\Omega_j$.
Across the interface $\Gamma:=\partial\Omega_1$, the continuity of the electric and magnetic fields leads to the condition
\begin{equation}\label{eq:int_cond}
u_1 = u_2, \quad  A_1 \nabla u_1 \cdot \nu = A_2 \nabla u_2 \cdot \nu, 
\end{equation}
where $\nu$ represents the unit outward normal along the interface $\Gamma$.

\subsection{Computational challenges}\label{sec:comput_challenge}
One unique feature of a hyperbolic medium is the strongly directional wave propagation inside a cone
with the half cone angle given by $\arctan\sqrt{-\dfrac{\Re\, \varepsilon_1}{\Re\, \varepsilon_2}}$.
This induces sharp transition of the solution across the cone boundary when wave starts decaying/growing exponentially.
The domain discretization methods such as finite element or finite difference are computationally very expensive to resolve such singular behaviors of the solution.
Here we propose a boundary element method for solving for the transmission problem \eqref{eq:trans_prob}.
Integral equation solvers have played an increasing role in the computational electromagnetics in the past several decades due to 
its powerfulness in solving large-scale problems by discretization over the boundaries of the objects only; 
see, for instance, the monographs \cite{habib-book, Chew1, Chew2, Colton-Kress, Hsiao-Wendeland, Nedelec} and the references therein. 
The application of the integral equation method for the hyperbolic media requires us to address several new computational challenges as described below.

\begin{figure}[!htbp]
\begin{center}
\includegraphics[height=6cm]{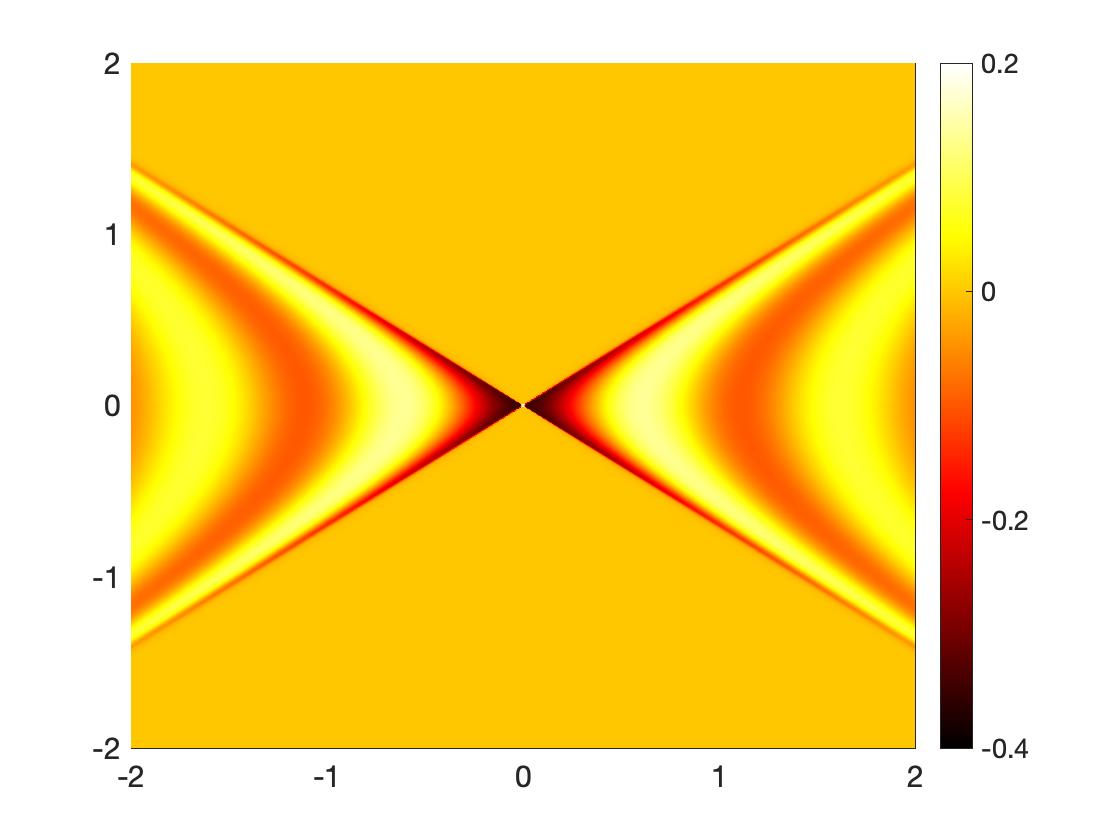}
\includegraphics[height=6cm]{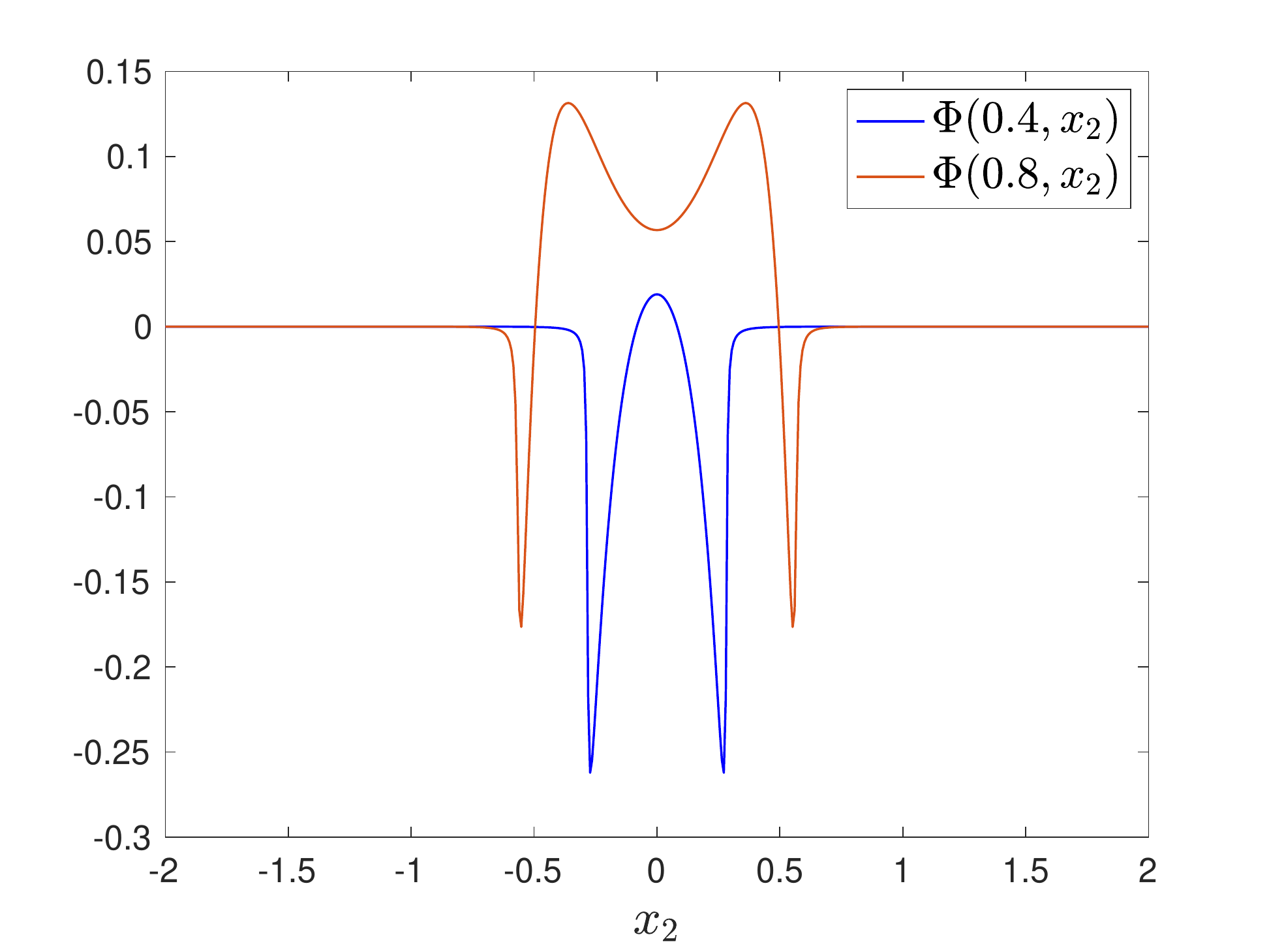}
\caption{Left: The real part of the fundamental solution $\Phi(x_1,x_2)$ in the hyperbolic medium with $\varepsilon_1=1+0.02i$ and $\varepsilon_2=-2+0.02i$. 
The source is located at the origin, and the wavenumber $k_0=2\pi$. Right: cross-sectional plot of  $\Phi(x_1,x_2)$ when $x_1=0.4$ and $0.8$.
$\Phi(x_1,x_2)$ attains sharp transitions near the boundary of the propagating cone $\mathcal C$.
 }\label{fig:green_fun}
\end{center}
\end{figure}

First, the fundamental solution
in the hyperbolic medium is strongly directional in the propagating cone 
$$ \mathcal C:=\{ x\,;\, x^T (\Re A)^{-1} x >0 \}, $$
and it decays exponentially across the cone boundary (see Figure \ref{fig:green_fun}).
As a result, when computing the stiff matrices in the boundary element method, one needs to apply adaptive numerical quadrature rules to
evaluate the integrals with the fundamental solution or its derivatives as kernels in order to achieve sufficient accuracy for discretization.
Second, the solution of the integral equation formulation along the interface $\Gamma$ attains sharp transitions
when wave front reaches the boundary, especially for hyperbolic media with small loss (see examples in Section \ref{sec:num_examples}).
Here we employ a two-level a posteriori error estimator and an adaptive mesh refinement procedure to resolve the singularity 
of the integral equation solution in an accurate and efficient manner. The theory and computation for the adaptive boundary element methods (BEM) are
mature in solving elliptic boundary value problems \cite{Feischl}. By using a posteriori error estimator,
the adaptive procedure chooses a sequence of meshes such that the numerical error decays in an optimal manner with increasing dimension of the approximation spaces.
There exist a variety of error estimators for elliptic boundary value problems, including residual type estimators, space enrichment type estimators,
averaging estimators, etc \cite{Carstensen1, Carstensen2, Carstensen3, Carstensen4, Faermann1, Faermann2, Ferraz-Leite1, Ferraz-Leite2, Wendland}.
The two-level a posteriori error estimator proposed here for the hyperbolic transmission problem belongs to the family of the space enrichment type estimators.
The principal idea is to use improved approximation of solutions $u_{h/2}$ and $\partial_{\tilde\nu}u_{h/2}$ obtained
over a uniform refined mesh with mesh size $h/2$ to replace the exact solution $u$ and $\partial_{\tilde\nu}u$ in the numerical error $\|u-u_h\|$ and $\|\partial_{\tilde\nu}u-\partial_{\tilde\nu} u_h\|$.
The D\"{o}rfler strategy is then applied to mark and refine the mesh where local errors $\|u_{h/2}-u_h\|+\|\partial_{\tilde\nu}u_{h/2}-\partial_{\tilde\nu} u_h\|$ are large.
Such an adaptive procedure allows for the reduction of the degree of freedom significantly while achieving desired accuracy for the solution, 
as demonstrated by the numerical examples in Section \ref{sec:num_examples}.
The goal of our work in this paper is to demonstrate the efficacy and accuracy of the adaptive algorithm for the two dimensional problems.
Its application in three dimensions will be investigated in the forthcoming work.

The rest of the paper is organized as follows. In Section \ref{sec:lp-BIE} we introduce layer potentials and derive their limits as the target points approach the boundary.
The limiting formulas recover the formulas in the isotropic medium when $\varepsilon_1=\varepsilon_2$.
The boundary integral equation for the transmission problem \eqref{eq:trans_prob} is then formulated in Section \ref{sec:lp-BIE}.
The adaptive Galerkin boundary element method is desribed in Section \ref{sec:adaptive-BEM}, where we introduce the adaptive numerical quadrature
and the two-level a posteriori error estimator. Several numerical examples are given in Section \ref{sec:num_examples} to illustrate the accuracy and efficiency of the adaptive algorithm.
The paper is concluded with brief remark about the proposed computational approach and the future work along this direction.

\section{Layer potentials and boundary integral equations for the transmission problem}\label{sec:lp-BIE}
\subsection{Layer potentials and integral operators}\label{sec:layer-potentials}
Here and henceforth,  for a hyperbolic material with $\varepsilon_\perp(x) \equiv \varepsilon_1$ and $\varepsilon_\parallel(x)\equiv\varepsilon_2$,
we let 
$$ \tilde r(x,y) = \sqrt{(x-y)^T A^{-1} (x-y)} $$
be the $A$-deformed distance between $x$ and $y$, where  $A$ is given in \eqref{eq:A} and
 the function $\sqrt{z}$ is understood as an analytic function defined in the domain $\mathbf{C}\backslash\{-it: t\geq 0\}$ such that
$\sqrt{z} = |z|^{\frac{1}{2}} e^{\frac{1}{2} i\arg z}$.  Note that when $\Im \varepsilon_j>0$ ($j=1,2$)$,  \tilde r(x,y)$ lies in the first quadrant of the complex plane.
Let $ \tilde \nu= A \nu$ be the $A$-deformed normal vector over the interface $\Gamma$.
Correspondingly, the derivative of a given function $\varphi$ along the direction $\tilde\nu$ is defined as
\begin{equation*}
\frac{\partial \varphi}{\partial  \tilde \nu} = \nabla \varphi \cdot \tilde \nu = \left[ \frac{1}{\varepsilon_1}\dfrac{\partial \varphi}{\partial x_1},    \frac{1}{\varepsilon_2}\dfrac{\partial \varphi}{\partial x_2} \right ]  \cdot \nu.
\end{equation*}

Let $\Phi(x,y) = \dfrac{i}{4} \sqrt{\varepsilon_1 \varepsilon_2}  H_0^{(1)}\left(k \tilde r(x,y)\right)$ be the fundamental solution,
which satisfies \eqref{eq:hyperbolic} when $x \neq y$ and is outgoing when $|x-y|\to\infty$.
Here $H_0^{(1)}(r)$ represents the zero order Hankel function of the first kind.
$\Omega$ is a bounded simply connected domain with the boundary $\Gamma$ of class $C^2$.
Given the density function $\varphi$ over $\Gamma$, the single and double layer potentials are defined by
\begin{equation}\label{eq:vw}
 v(x) = \int_\Gamma \Phi(x,y) \varphi(y) \, ds_y \quad \mbox{and} \quad
 w(x) = \int_\Gamma  \frac{\partial \Phi(x,y)}{\partial \tilde\nu(y)} \varphi(y) \, ds_y
\quad \mbox{for} \quad x\notin\Gamma.
\end{equation}
It is well-known that the single layer potential $v(x)$ is continuous throughout $\mathbb{R}^2$.
In what follows, we derive the limits of the double layer potential and the derivatives of two layer potentials as $x$ approaches $\Gamma$. 
The limiting formulas recover the classical limiting formulas when $\varepsilon_1=\varepsilon_2=1$ \cite{Kress}.

\begin{lem}\label{lem:w}
The double-layer potential $w(x)$ with the continuous density $\varphi$  can be continuously extended 
from $\Omega$ to $\bar\Omega$ and $\mathbb{R}^2\backslash\bar\Omega$ to $\mathbb{R}^2\backslash\Omega$ respectively with the limit
\begin{equation}\label{eq:w_limit}
 w_\pm (x) = \int_\Gamma  \frac{\partial \Phi(x,y)}{\partial \tilde\nu(y)} \varphi(y) \, ds_y \pm \frac{1}{2}  \varphi(x)  \quad \mbox{for} \; x\in \Gamma,
\end{equation}
where 
$$w_\pm(x) = \lim_{h>0, h\to0} w(x\pm h\nu(x)).$$
\end{lem}
\noindent\textbf{Proof.} Let $\Phi_0(x,y) = -\frac{1}{2\pi} \sqrt{\varepsilon_1 \varepsilon_2}  \ln (\tilde r(x,y))$ be the fundamental solution of \eqref{eq:hyperbolic} when $k_0=0$, and $w_0$
be the corresponding double layer potential:
$$  w_0(x) = \int_\Gamma  \frac{\partial \Phi_0(x,y)}{\partial \tilde\nu(y)} \varphi(y) \, ds_y. $$
Note that the difference of two double layer potentials $w(x)$ and $w_0(x)$ is continuous in $\mathbb{R}^2$, thus it suffices to verify \eqref{eq:w_limit} for $w_0(x)$. 
The proof can be further reduced to the special case when the density function $\varphi\equiv 1$, this is because for an arbitrary density function $\varphi$,
one can write the double layer potential as
$$  w_0(x) = \varphi(x)  \int_\Gamma  \frac{\partial \Phi_0(x,y)}{\partial \tilde\nu(y)}  \, ds_y + \int_\Gamma  \frac{\partial \Phi_0(x,y)}{\partial \tilde\nu(y)} (\varphi(y) - \varphi(x)) \, ds_y,  $$
and the latter is continuous  throughout $\mathbb{R}^2$ when $\varphi$ is continuous. Next we verify the assertion by assuming that $\varphi(x)\equiv1$ and showing that
\begin{equation}\label{eq:w0}
w_0(x)=\left\{
\begin{array}{ccc}
0, & x\in \mathbb{R}^2\backslash\bar\Omega, \\
-\frac{1}{2}, & x\in \Gamma, \\
-1, & x\in \Omega.
\end{array}\right.
\end{equation}

When $x\in \mathbb{R}^2\backslash\bar\Omega$, noting that $\Phi_0(x,y)$ solves \eqref{eq:hyperbolic} with $k_0=0$ in $\Omega$,
it is obvious $w_0(x)=0$ by applying the Green's formula. Now if $x\in\Gamma$, let $B_\delta(x)$ be the small disk with radius $\delta$ centered at $x$
and $\partial B_\delta(x)$ be its boundary. It follows from the Green's formula that
\begin{equation*}
 \int_{\Gamma}  \frac{\partial \Phi_0(x,y)}{\partial \tilde\nu(y)}  \, ds_y = \lim_{\delta\to0} \int_{\partial B_\delta(x) \cap \Omega}  \frac{\partial \Phi_0(x,y)}{\partial \tilde\nu(y)}  \, ds_y, 
 \end{equation*}
where $\nu$ denotes the unit normal exterior to $B_\delta(x)$. A direct calculation yields 
\begin{equation*}
\int_{\partial B_\delta(x) \cap \Omega}  \frac{\partial \Phi_0(x,y)}{\partial \tilde\nu(y)}  \, ds_y = - \frac{1}{2\pi}  \sqrt{\varepsilon_1 \varepsilon_2} \int_{\theta_1(\delta)}^{\theta_2(\delta)} \frac{r^2}{\tilde r^2}  \, d\theta
= - \frac{1}{2\pi} \int_{\theta_1(\delta)}^{\theta_2(\delta)} \frac{ \sqrt{\varepsilon_1 \varepsilon_2}  }{  \varepsilon_1 \cos^2\theta + \varepsilon_2 \sin^2\theta}  \, d\theta
\end{equation*}
in the polar coordinate,  where  $r^2 = (x_1-y_1)^2 +  (x_2-y_2)^2$ and $\tilde r ^2= \varepsilon_1 (x_1-y_1)^2 + \varepsilon_2 (x_2-y_2)^2 $.
Hence, there holds
\begin{equation*}
 \int_{\Gamma}  \frac{\partial \Phi_0(x,y)}{\partial \tilde\nu(y)}  \, ds_y = \lim_{\delta\to0}    - \frac{1}{2\pi} \sqrt{\frac{\varepsilon_1}{\varepsilon_2}} \int_{\theta_1(\delta)}^{\theta_2(\delta)} \frac{  \sec^2\theta }{ \tan^2\theta +\varepsilon_1/\varepsilon_2 }  \, d\theta.
 \end{equation*}
By evaluating the above integral explicitly and noting that $\lim_{\delta\to0} (\theta_2(\delta) - \theta_1(\delta)) = \pi $, we obtain that $w_0(x)=-\frac{1}{2}$.
A parallel calculation leads to $w_0(x)=-1$ for $x\in\Omega$, for which there holds $\theta_2(\delta) - \theta_1(\delta) = 2\pi$.  \qed

\begin{prop}\label{lem:dv}
The derivative of the single-layer potential $v(x)$ with the continuous density $\varphi$  can be continuously extended 
from $\Omega$ to $\bar\Omega$ and $\mathbb{R}^2\backslash\bar\Omega$ to $\mathbb{R}^2\backslash\Omega$ respectively with the limit
\begin{equation}\label{eq:dv_limit}
\frac{\partial v_\pm(x)}{\partial \tilde\nu} = \int_\Gamma  \frac{\partial \Phi(x,y)}{\partial \tilde\nu(x)} \varphi(y) \, ds_y \mp \frac{1}{2}  \varphi(x)  \quad \mbox{for} \; x\in \Gamma,
\end{equation}
where $$\frac{\partial v_\pm(x)}{\partial \tilde\nu} = \lim_{h>0,h\to0}  \nabla v(x\pm h\nu(x)) \cdot \tilde\nu(x).$$
\end{prop}
\noindent\textbf{Proof.}  This can be observed from the formula
$$ \nabla v(x\pm h\nu(x)) \cdot \tilde\nu(x) = -w(x\pm h\nu(x)) + \int_\Gamma   \nabla_y\Phi(x\pm h\nu(x),y) \cdot (\tilde\nu(y) - \tilde\nu(x)  ) \varphi(y) \, ds_y,  $$
wherein we have used the relation $ \nabla_x \Phi(x,y)=- \nabla_y \Phi(x,y)$.
Note that the second term is continuous in $\mathbb{R}^2$, thus an application of Lemma \ref{lem:w} leads to \eqref{eq:dv_limit}.      \qed

\medskip

\begin{lem}\label{lem:gradv}
The gradient of the single layer potential $v(x)$ with the density $\varphi \in C^1(\Gamma)$ can be continuously extended
from $\Omega$ to $\bar\Omega$ and $\mathbb{R}^2\backslash\bar\Omega$ to $\mathbb{R}^2\backslash\Omega$ respectively with the limit
\begin{equation}\label{eq:gradv_limit}
 \nabla v_\pm (x) = \int_\Gamma  \nabla_x\Phi(x,y) \varphi(y) \, ds_y \mp \frac{1}{2}  \varphi(x) \, e_{\nu(x)} \quad \mbox{for} \; x\in \Gamma.
\end{equation}
In the above, $\displaystyle{\nabla  v_\pm(x) = \lim_{h>0, h\to0} \nabla v(x\pm h\nu(x))}$, and the vector $e_{\nu(x)}$ is given by
$e_{\nu(x)}= \dfrac{\nu(x)}{\nu(x) \cdot \tilde\nu(x)} $.
\end{lem}

\noindent\textbf{Proof.}  For $y\in\Gamma$, let $\nu(y)=[\nu_1(y), \nu_2(y)]^T$ and $\tau(y)=[-\nu_2(y), \nu_1(y)]^T$ be the normal and tangential vector respectively.
For $x\neq y$, one can decompose $\nabla_y \Phi(x,y)$ as
\begin{equation}\label{eq:Phi_decomp}
 \nabla_y \Phi(x,y) = \frac{\partial \Phi(x,y)}{\partial \tilde\nu(y)} \, e_{\nu(y)} +  \frac{\partial \Phi(x,y)}{\partial \tau(y)}  \, e_{\tau(y)}, 
\end{equation}
where
\begin{equation*}
e_{\nu(y)}= \frac{1}{\nu(y) \cdot \tilde\nu(y)} \left[
\begin{array}{ccc}
\nu_1(y)  \\
\nu_2(y)
\end{array} \right]
\quad \mbox{and} \quad
e_{\tau(y)}= \frac{1}{\nu(y) \cdot \tilde\nu(y)} \left[
\begin{array}{ccc}
-\nu_2(y)/\varepsilon_2  \\
\nu_1(y)/\varepsilon_1
\end{array} \right].
\end{equation*}

For fixed $z\in\Gamma$, set $x=z+h\nu(x)$ with $0<h\ll1$. Note that $ \nabla_x \Phi(x,y)=- \nabla_y \Phi(x,y)$, and in light of the decomposition \eqref{eq:Phi_decomp},
we have
\begin{equation*}
\int_\Gamma \nabla_x \Phi(x,y) \varphi(y) \, ds_y = - \int_\Gamma  \frac{\partial \Phi(x,y)}{\partial \tilde\nu(y)} \varphi(y)e_{\nu(y)}   \, ds_y
- \int_\Gamma  \frac{\partial \Phi(x,y)}{\partial \tau(y)} \varphi(y)e_{\tau(y)}   \, ds_y
\end{equation*}
By letting $h\to0$, there holds
\begin{eqnarray*}
 \lim_{h\to0} \int_\Gamma  \frac{\partial \Phi(x,y)}{\partial \tilde\nu(y)} \varphi(y)e_{\nu(y)}   \, ds_y &=& \int_\Gamma  \frac{\partial \Phi(z,y)}{\partial \tilde\nu(y)} \varphi(y)e_{\nu(y)}   \, ds_y + \frac{1}{2}  \varphi(z) \cdot e_{\nu(z)}, \\
\lim_{h\to0}  \int_\Gamma  \frac{\partial \Phi(x,y)}{\partial \tau(y)} \varphi(y)e_{\tau(y)} \, ds_y &=& - \lim_{h\to0} \int_\Gamma  \Phi(x,y) \dfrac{d (\varphi(y)e_{\tau(y)})}{d\tau(y)} \, ds_y \\
&=& -\int_\Gamma  \Phi(z,y) \dfrac{d(\varphi(y)e_{\tau(y)})}{d\tau(y)} \, ds_y  \\
&=& \int_\Gamma  \frac{\partial \Phi(z,y)}{\partial \tau(y)} \varphi(y)e_{\tau(y)} \, ds_y.
\end{eqnarray*}
Therefore,
\begin{eqnarray*}
\lim_{h\to0} \int_\Gamma \nabla_x \Phi(x,y) \varphi(y) \, ds_y &=& -\int_\Gamma  \left(\frac{\partial \Phi(z,y)}{\partial \tilde\nu(y)} e_{\nu(y)}  +   \frac{\partial \Phi(z,y)}{\partial \tau(y)} e_{\tau(y)} \right) \varphi(y) \, ds_y - \frac{1}{2}  \varphi(z) \cdot e_{\nu(z)} \\
&=& \int_\Gamma  \nabla_z\Phi(z,y) \varphi(y) \, ds_y - \frac{1}{2}  \varphi(z) \cdot e_{\nu(z)},
\end{eqnarray*}
where we use \eqref{eq:Phi_decomp} and the relation $ \nabla_z \Phi(z,y)=- \nabla_y \Phi(z,y)$ again. The proof for $x=z-h\nu(x)$ is parallel.

\qed

\begin{lem}\label{lem:der_dbl}
Let $w(x)$ be the double layer potential defined in \eqref{eq:vw} and $\varphi \in C^1(\Gamma)$, then there holds
\begin{equation}\label{eq:dw_limit}
\frac{\partial w_\pm(x)}{\partial \tilde\nu} = \frac{1}{\varepsilon_1\varepsilon_2}\dfrac{d}{ds_x} \int_\Gamma  \Phi(x,y) \frac{d\varphi(y)}{ds} \, ds_y  
+ k_0^2 \int_\Gamma \Phi(x,y) \varphi(y) (\nu(x)  \cdot \tilde \nu(y)) \, ds_y   \quad \mbox{for} \; x\in \Gamma,
\end{equation}
where  $\displaystyle{ \frac{\partial w_\pm(x)}{\partial \tilde\nu} = \lim_{h>0, h\to0}  \nabla w(x\pm h\nu(x)) \cdot \tilde\nu(x)} $,
and $\frac{d}{ds}$ denotes the derivative with respect to the arc length.
\end{lem}
\noindent\textbf{Proof.}  For $y\in\Gamma$, let $\tau(y)^\perp = [ \tau_2(y), -\tau_1(y) ]^T = \nu(y)$.
If $x\not\in\Gamma$, using the relation
$$  \frac{\partial \Phi(x,y)}{\partial \tilde\nu(y)} =  \nabla_y\Phi(x,y) \cdot \tilde\nu(y) = -\nabla_x\Phi(x,y) \cdot \tilde\nu(y)
 =  - \left[ \frac{1}{\varepsilon_1}\frac{\partial\Phi(x,y)}{\partial x_1},  \;  \frac{1}{\varepsilon_2}\frac{\partial\Phi(x,y)}{\partial x_2} \right ] \cdot  \tau(y)^\perp, $$
we have
\begin{eqnarray*}
\nabla_x \frac{\partial \Phi(x,y)}{\partial \tilde\nu(y)} &=& - \nabla_x \left(  \frac{\tau_2(y) }{\varepsilon_1}  \frac{\partial\Phi(x,y)}{\partial x_1} -  \frac{\tau_1(y) }{\varepsilon_2}  \frac{\partial\Phi(x,y)}{\partial x_2} \right) \\
&=&  \left[ \frac{\tau_1(y) }{\varepsilon_2}  \frac{\partial^2 \Phi(x,y)}{\partial x_1\partial x_2}  -  \frac{\tau_2(y) }{\varepsilon_1}  \frac{\partial^2\Phi(x,y)}{\partial x_1^2}, \; \frac{\tau_1(y) }{\varepsilon_2}  \frac{\partial^2\Phi(x,y)}{\partial x_2^2}   -  \frac{\tau_2(y)}{\varepsilon_1} \frac{\partial^2\Phi(x,y)}{\partial x_1\partial x_2}  \right]^T.
\end{eqnarray*}
By using the equation \eqref{eq:hyperbolic} and the relation $\nabla_x\Phi(x,y)=-\nabla_y\Phi(x,y)$, it follows that
\begin{eqnarray*}
\nabla_x \frac{\partial \Phi(x,y)}{\partial \tilde\nu(y)}
&=& \left[ \frac{\tau_1(y) }{\varepsilon_2}  \frac{\partial^2 \Phi(x,y)}{\partial x_1\partial x_2}  +  \frac{\tau_2(y) }{\varepsilon_2}  \frac{\partial^2\Phi(x,y)}{\partial x_2^2}, \; 
               -\frac{\tau_1(y) }{\varepsilon_1}  \frac{\partial^2\Phi(x,y)}{\partial x_1^2}   -  \frac{\tau_2(y)}{\varepsilon_1} \frac{\partial^2\Phi(x,y)}{\partial x_1\partial x_2}  \right]^T \\
               && + k_0^2 \Phi(x,y)  \tau(y)^\perp \\
&=& \left[ -\frac{1}{\varepsilon_2}  \frac{\partial}{\partial x_2} \frac{\partial \Phi(x,y) }{\partial \tau_y}, \;  \frac{1}{\varepsilon_1}  \frac{\partial}{\partial x_1} \frac{\partial \Phi(x,y) }{\partial \tau_y}\right]^T
+k_0^2 \Phi(x,y)  \nu(y).
\end{eqnarray*}
Hence, applying the integration by parts leads to
\begin{equation}\label{eq:der_db2}
A  \nabla_x  \int_\Gamma \frac{\partial \Phi(x,y)}{\partial \tilde\nu(y)} \,  \varphi(y) ds_y 
= \frac{1}{\varepsilon_1\varepsilon_2} \left(\nabla_x \int_\Gamma  \Phi(x,y) \frac{d\varphi(y)}{ds} \, ds_y \right)^\perp 
 + k_0^2 \int_\Gamma \Phi(x,y) \varphi(y)  \tilde\nu(y) \, ds_y.
\end{equation}

Now for fixed $z\in\Gamma$, we set $x=z+h\nu(z)$. By virtue of \eqref{eq:der_db2} and Lemma \ref{lem:gradv}, it follows that
\begin{eqnarray*}
&& \lim_{h\to0}  \nabla w(z + h\nu(z)) \cdot \tilde\nu(z) \\
&=&  \lim_{h\to0}  \frac{1}{\varepsilon_1\varepsilon_2} \left(\nabla \int_\Gamma  \Phi(z+h\nu(z), y) \frac{d\varphi(y)}{ds} \, ds_y \right)^\perp \cdot \nu(z)
 + k_0^2 \int_\Gamma \Phi(z,y) \varphi(y)  (\nu(z)\cdot\tilde\nu(y)) \, ds_y \\
 &=&  \frac{1}{\varepsilon_1\varepsilon_2} \left(\nabla_z \int_\Gamma  \Phi(z, y) \frac{d\varphi(y)}{ds} \, ds_y - \frac{1}{2}  \varphi(z) \, e_{\nu(z)} \right)^\perp \cdot \nu(z)
 + k_0^2 \int_\Gamma \Phi(z,y) \varphi(y)  (\nu(z)\cdot\tilde\nu(y)) \, ds_y.
 \end{eqnarray*}
A straightforward calculation leads to 
$$ \frac{\partial w_+(z)}{\partial \tilde\nu} = \frac{1}{\varepsilon_1\varepsilon_2}\dfrac{d}{ds_z} \int_\Gamma  \Phi(z,y) \frac{d\varphi(y)}{ds} \, ds_y  
+ k_0^2 \int_\Gamma \Phi(z,y) \varphi(y) (\nu(z)  \cdot \tilde \nu(y)) \, ds_y. $$
Similarly, by setting $x=z-h\nu(z)$, one can obtain the same formula for $\frac{\partial w_-(z)}{\partial \tilde\nu} $.   \qed

\subsection{Boundary integral equation formulation for the transmission problem}\label{sec:BIE}
Let $\Phi_j(x,y) = \dfrac{i}{4} \sqrt{\varepsilon_1^{(j)} \varepsilon_2^{(j)}}  H_0^{(1)}\left(k \tilde r_j(x,y)\right)$ be the fundamental solution in the domain $\Omega_j$,
where $\tilde r_j = (x-y)^T A_j^{-1} (x-y) $.  Let  $ \tilde \nu_j= A_j \nu$ be the $A_j$-deformed normal vector over the interface $\Gamma$.
We define the integral operators $S_j$, $K_j$, $K'_j$ and $N_j$ for $x\in \Gamma$ as follows:
\begin{align}
 [S_j \varphi] (x) & = \int_\Gamma \Phi_j(x,y) \varphi(y) \, ds_y,  \label{eq:Sj} \\
 [K_j \varphi] (x) & = \int_\Gamma \frac{\partial \Phi_j(x,y)}{\partial \tilde\nu_j(y)} \varphi(y) \, ds_y,  \label{eq:Kj}  \\
 [K'_j \varphi] (x) & = \int_\Gamma \frac{\partial \Phi_j(x,y)}{\partial \tilde\nu_j(x)} \varphi(y) \, ds_y,  \label{eq:K'j} \\
 [N_j \varphi] (x) & = \frac{1}{\varepsilon_1^{(j)} \varepsilon_2^{(j)}}\dfrac{d}{ds_x} S_j \left(\frac{d\varphi}{ds}\right)  + k_0^2  \nu(x) \cdot S_j (\varphi\tilde \nu_j). \label{eq:Nj} \\
& = \frac{1}{\varepsilon_1^{(j)} \varepsilon_2^{(j)}}\dfrac{d}{ds_x} \int_\Gamma  \Phi_j(x,y) \frac{d\varphi(y)}{ds} \, ds_y
+ k_0^2 \int_\Gamma \Phi_j(x,y) \varphi(y) (\nu(x)  \cdot \tilde \nu_j(y)) \, ds_y. \nonumber
\end{align}
Note that $\Phi_j(x,y)$ in the lossy hyperbolic medium attains the same singularity as the fundamental solution with $\varepsilon_1=\varepsilon_2=1$,  
thus from the standard theory of the boundary integral operators (cf. \cite{Kress}),
we have the following lemma for the above integral operators.
\begin{lem}\label{lem:opts}
The operators $S_j: H^{-1/2}(\Gamma) \to H^{1/2}(\Gamma)$,  $K_j: H^{1/2}(\Gamma) \to H^{1/2}(\Gamma)$, 
$K_j': H^{-1/2}(\Gamma) \to H^{-1/2}(\Gamma)$, and $N_j: H^{1/2}(\Gamma) \to H^{-1/2}(\Gamma)$ are bounded.
\end{lem}

Let
$$ S=S_1-S_2,  \quad K=K_1-K_2, \quad K'=K_1' -K_2', \quad N =N_1' -N_2' $$
be the difference of two integral operators with associated kernels.
The volume integral operators $P_j$ and $Q_j$ are defined as
\begin{align*}
[P_j \varphi] (x) =  \int_{\Omega_j} \Phi_j(x,y) \varphi(y) \, dy, \quad \mbox{for} \; x\in \Gamma, \\
[Q_j \varphi] (x) = \int_{\Omega_j} \frac{\partial \Phi_j(x,y)}{\partial \nu(x)} \varphi(y) \, dy, \quad \mbox{for} \; x\in \Gamma.
\end{align*}

Applying the Green's formula in $\Omega_1$ and using the formula \eqref{eq:w0}, we obtain for $x\in\Omega_1$ that
\begin{equation}\label{eq:u_omega1}
u_1(x) = \int_\Gamma  \Phi_1(x,y)  \frac{\partial u_1(y)}{\partial \tilde\nu_1(y)}  -     \frac{\partial \Phi_1(x,y)}{\partial \tilde\nu_1(y)} u_1(y)  \, ds_y - \int_{\Omega_1}  \Phi_1(x,y)  f_1(y) \, dy.  \\
\end{equation}
Similarly,  for $x\in\Omega_2$ there holds
\begin{equation}\label{eq:u_omega2}
u_2(x) = \int_\Gamma  \frac{\partial \Phi_2(x,y)}{\partial \tilde\nu_2(y)} u_2(y) - \Phi_2(x,y)  \frac{\partial u_2(y)}{\partial \tilde\nu_2(y)}  \, ds_y - \int_{\Omega_2}  \Phi_2(x,y)  f_2(y) \, dy .
\end{equation}
Recall that $f_1(x)$ and $f_2(x)$ are localized sources with compact support, the volume integrals above only need to be evaluated over their support regions.
By taking the limit of  \eqref{eq:u_omega1} and \eqref{eq:u_omega2} when $x$ approaches the interface $\Gamma$ and applying Lemma \ref{lem:w},
we achieve the integral equations on $\Gamma$:
\begin{eqnarray}
\frac{1}{2}u_1&=& S_1  \left(\frac{\partial u_1}{\partial \tilde\nu_1} \right) - K_1 u_1  - P_1 f_1;  \label{eq:BIE1} \\
\frac{1}{2}u_2 &=& K_2 u_2 - S_2  \left(\frac{\partial u_2}{\partial \tilde\nu_2} \right)   -  P_2 f_2. \label{eq:BIE2}
\end{eqnarray}
For $x\in\Gamma$, evaluating $\nabla u_j(x\pm h\nu(x)) \cdot \tilde\nu(x)$ in \eqref{eq:u_omega1} and \eqref{eq:u_omega2} respectively, and taking the limit when $h\to0$
yields
\begin{eqnarray}
\frac{1}{2}\frac{\partial u_1}{\partial\tilde\nu} &=& K_1'  \left(\frac{\partial u_1}{\partial \tilde\nu_1} \right) - N_1 u_1  - Q_1 f_1;   \label{eq:BIE3} \\
\frac{1}{2}\frac{\partial u_2}{\partial\tilde\nu}&=& N_2 u_2 - K_2'  \left(\frac{\partial u_2}{\partial \tilde\nu_2} \right) -  Q_2 f_2, \label{eq:BIE4}
\end{eqnarray}
where we have used Lemmas \ref{lem:gradv} and \ref{lem:der_dbl}. 

By taking the sum \eqref{eq:BIE1} + \eqref{eq:BIE2} and  \eqref{eq:BIE3} + \eqref{eq:BIE4} respectively, and applying the continuity condition \eqref{eq:int_cond} for the wave field across the interface,
we obtain the following system of integral equations:
\begin{equation}\label{eq:BIE_tran}
\left[
\begin{array}{cc}
N &  I-K' \\
I+K  & -S 
\end{array} \right]
\left[
\begin{array}{cc}
\phi^{(1)}  \\
\phi^{(2)} 
\end{array} \right]
=
\left[
\begin{array}{cc}
 g^{(1)} \\
 g^{(2)}
\end{array} \right].
\end{equation}
In the above, $\phi^{(1)} = u_1 = u_2$, $\phi^{(2)} = \frac{\partial u_1}{\partial \tilde\nu_1} = \frac{\partial u_2}{\partial \tilde\nu_2}$, and 
$$ g^{(1)} = -(Q_1f_1+Q_2f_2), \quad g^{(2)} = -(P_1f_1+P_2f_2). $$
We point out that integral equations in the form of \eqref{eq:BIE_tran} have been widely used in studies of acoustic, electromagnetic, and elastic transmission problems
with isotropic media; see, for instance, \cite{Bu-Lin-Reitich, Costabel-Stephan, Kleinman-Martin} and the references therein.

For brevity of notation, we let 
\begin{equation}\label{eq:optT}
\mathbf{T} = \left[
\begin{array}{cc}
N &  I-K' \\
I+K  & -S 
\end{array} \right],
\quad
\boldsymbol{\phi} = 
\left[
\begin{array}{cc}
\phi^{(1)}  \\
\phi^{(2)} 
\end{array} \right],
\quad
\boldsymbol{g} =
\left[
\begin{array}{cc}
 g^{(1)} \\
 g^{(2)}
\end{array} \right].
\end{equation}
In view of Lemma \ref{lem:opts}, the operator $\mathbf{T}$ is bounded from $H^{1/2}(\Gamma) \times H^{-1/2}(\Gamma) \to H^{-1/2}(\Gamma) \times H^{1/2}(\Gamma) $.
Let $\mathcal{H}=H^{1/2}(\Gamma) \times H^{-1/2}(\Gamma)$ be the Hilbert space equipped with the inner product
$$ 
 \langle \boldsymbol{\phi}, \boldsymbol{\psi} \rangle = \int_\Gamma \phi^{(1)}\bar\psi^{(1)} + \phi^{(2)}\bar\psi^{(2)} \, ds,
$$
and $\mathcal{H}'=H^{-1/2}(\Gamma) \times H^{1/2}(\Gamma)$ be the dual space of $\mathcal{H}$. Then the weak
formulation for the integral equation equation \eqref{eq:BIE_tran} reads as finding $\boldsymbol{\phi}\in \mathcal{H}$ such that 
\begin{equation}\label{eq:weak_form}
a(\boldsymbol{\phi},\boldsymbol{\psi}) = \langle \boldsymbol{g}, \boldsymbol{\psi} \rangle \quad \forall \boldsymbol{\psi} \in \mathcal{H},
\end{equation}
in which the sesquilinear form 
\begin{equation*}%\label{eq:sesqlinear_form}
a(\boldsymbol{\phi},\boldsymbol{\psi}) = \langle \mathbf{T}\boldsymbol{\phi}, \boldsymbol{\psi} \rangle 
= \langle N\phi^{(1)}, \psi^{(1)} \rangle + \langle (I-K')\phi^{(2)}, \psi^{(1)} \rangle + \langle (I+K)\phi^{(1)}, \psi^{(2)} \rangle - \langle S\phi^{(2)}, \psi^{(2)} \rangle.
\end{equation*}

\section{Adaptive Galerkin boundary element method}\label{sec:adaptive-BEM}
\subsection{Galerkin boundary element method}
Let $\Gamma_h=\{ E_1, E_2, \cdots, E_M \}$ be a mesh assigned on the interface $\Gamma$ with the mesh size $h=\max_{1\le m\le M}|E_m|$,
and $\mathcal{H}_h:=\mathcal{U}_h \times \mathcal{V}_h$ be the corresponding finite-dimensional space defined over $\Gamma_h$.
Here $\mathcal{U}_h$ and $\mathcal{V}_h$ are the  finite-dimensional approximation of the space $H^{1/2}(\Gamma)$ and $H^{-1/2}(\Gamma)$
respectively. The Galerkin boundary element method is to find $\boldsymbol{\phi}_h\in\mathcal{H}_h$ such that 
\begin{equation}\label{eq:galerkin_form}
a(\boldsymbol{\phi}_h,\boldsymbol{\psi}_h) = \langle \boldsymbol{g}, \boldsymbol{\psi}_h \rangle \quad \forall \boldsymbol{\psi}_h \in \mathcal{H}_h.
\end{equation}

Let $\left\{ \phi_m^{(1)} \right\}_{m=1}^{M_1}$  and $\left\{ \phi_m^{(2)} \right\}_{m=1}^{M_2}$ be the basis of the space $\mathcal{U}_h$ and $\mathcal{V}_h$ respectively.
By representing the solution $\displaystyle \phi_h^{(1)}(x)=\sum_{m=1}^{M_1} c_m^{(1)} \phi_m^{(1)} (x) $ and $\displaystyle \phi_h^{(2)}(x)=\sum_{m=1}^{M_2} c_m^{(2)} \phi_m^{(2)}(x) $, 
the Galerkin  approximation \eqref{eq:galerkin_form} leads to the following linear system:
\begin{equation}\label{eq:linear_system}
\left[
\begin{array}{cc}
N_h &  I_h'-K_h' \\
I_h+K_h  & -S_h 
\end{array} \right]
\left[
\begin{array}{cc}
\mathbf{c}^{(1)}  \\
\mathbf{c}^{(2)} 
\end{array} \right]
=
\left[
\begin{array}{cc}
 g_h^{(1)} \\
 g_h^{(2)}
\end{array} \right].
\end{equation}
In the above, $\mathbf{c}^{(1)}$ and $\mathbf{c}^{(2)}$ represent the unknown vectors that take the following form:
$$
\mathbf{c}^{(1)} = \left[ c_1^{(1)}, c_2^{(1)}, \cdots, c_{M_1}^{(1)} \right], \quad
\mathbf{c}^{(2)} = \left[ c_1^{(2)}, c_2^{(2)}, \cdots, c_{M_2}^{(2)} \right].
$$
The $(m,n)$-th entry of the matrices $N_h$,  $K_h$,  $S_h$, and $I_h$ are given by
\begin{eqnarray}
N_h(m,n) &=& \left\langle N\phi_n^{(1)}, \phi_m^{(1)}  \right\rangle, \quad K_h(m,n) = \left\langle K\phi_n^{(1)}, \phi_m^{(2)}  \right\rangle,  \label{eq:NK_h}  \\
S_h(m,n) &=& \left\langle S\phi_n^{(2)}, \phi_m^{(2)}  \right\rangle, \quad I_h(m,n) = \left\langle \phi_n^{(1)}, \phi_m^{(2)}  \right\rangle.    \label{eq:SI_h} 
\end{eqnarray}
$I_h'$ and $K_h'$ are the transposes of $I_h$ and $K_h$ respectively.
In light of \eqref{eq:Nj}, we evaluate $N_h(m,n)$ using the formula
\begin{equation}\label{eq:Nj_evaluation}
 \left\langle N_j\phi_n^{(1)}, \phi_m^{(1)} \right\rangle = -\frac{1}{\varepsilon_1^{(j)} \varepsilon_2^{(j)}} \left\langle S_j \frac{d\phi_n^{(1)}}{ds}, \frac{d\phi_m^{(1)}}{ds} \right\rangle
+  k_0^2  \left\langle S_j (\phi_n^{(1)} \tilde \nu_j)  \cdot  \nu,  \phi_m^{(1)}   \right\rangle, \quad j=1,2.
\end{equation}
The $m$-th element for the vectors $g_h^{(1)}$ and $g_h^{(2)}$ are given by
$$  g_h^{(1)}(m) =  \left\langle  g^{(1)}, \phi_m^{(1)}  \right\rangle, \quad  g_h^{(2)}(m) =  \left\langle  g^{(2)}, \phi_m^{(2)}  \right\rangle. $$

\subsection{Adaptive numerical integration}\label{sec:adapt_num_int}
The entries of the local stiff matrices \eqref{eq:NK_h} - \eqref{eq:Nj_evaluation}  boil down to the evaluation of integral in the form of
\begin{equation*}
e_{mn} = \int_{E_m} \int_{E_n} \Theta(x,y) \phi_n(y) \, \phi_m(x) \ ds_y ds_x,
\end{equation*}
in which $\Theta(x,y)$ represents the kernel $\Phi_j(x,y)$ or $\dfrac{\partial \Phi_j(x,y)}{\partial \tilde\nu_j(y)}$, and $\phi_m(x)$ represents the basis functions
in $\mathcal{U}_h$ or $\mathcal{V}_h$. 
As pointed out in Section \ref{sec:comput_challenge}, 
the fundamental solution in the hyperbolic medium is strongly directional in the propagating cone $\mathcal C=\{ x\,;\, x^T (\Re A)^{-1} x >0 \}$,
and it attains sharp transition near the cone boundary $\partial \mathcal C$. To capture the variation of the kernels accurately, we employ the adaptive numerical quadrature
to compute $e_{mn}$.  In more details, for a given small real number $\tau>0$, let us introduce the domain
\begin{equation}\label{eq:Omega_tau}
 \Omega_\tau=\{ x\,;\,  dist(x, \partial \mathcal C) < \tau \}
\end{equation}
that includes the cone boundary. $\tau$ is chosen so that $\Omega_\tau$ contains the region where $\Theta(x,y)$ attains very large derivatives.
Let $E_{mn} =\{ x-y\,;\, x\in E_m, y\in E_n \}$ be the set of relative locations between the source and target points when evaluating $e_{mn}$.
If $E_{mn}\cap\Omega_\tau=\emptyset$, then $\Theta(x,y)$ changes smoothly over the domain $(x,y)\in E_m \times E_n$,
thus there is no need for adaptivity and one can still apply fast algorithms, such as the
fast multipole methods, to evaluate $e_{mn}$ for fixed $m$ and all $n$ satisfying $E_{mn}\cap\Omega_\tau=\emptyset$ in an efficient manner \cite{Chew1, Greengard-Rokhlin}. 
Otherwise, if $E_{mn}\cap\Omega_\tau\neq\emptyset$, we compute $e_{mn}$ adaptively to resolve the kernels accurately as described below.

Note that when $m=n$, the kernel $\Theta(x,y)$ is weakly singular and the singular part of $e_{mn}$ can be evaluated analytically.
Hence we only need to consider the case when $\Theta(x,y)$ is nonsingular. By a change of variable,  $e_{mn}$ is expressed as
\begin{equation}\label{eq:emn}
e_{mn} = \int_{-1}^1 \int_{-1}^1 \Theta \big(x(t),y(s)\big) \phi_n\big(y(s)\big) \ ds \, \phi_m\big(x(t)\big)  dt
\end{equation}
in the parameter space. The integral in \eqref{eq:emn} is computed via the adaptive Lobatto quadrature rule \cite{Gander-Gautschi}. 
Let $\left\{ R_j^{(\ell)} \right\}_{j=1}^J$ be a decomposition of the whole region at level $\ell$ with small rectangles $R_j^{(\ell)}$. 
Starting from $\ell=0$ and $J=1$ with $R_1^{(0)}=[-1,1]\times[-1,1]$,
the adaptive algorithm computes the integral recursively over the region $R_j^{(\ell)}$ by first dividing $R_j^{(\ell)}$ in half along $t$ and $s$ coordinate axis respectively to produce four new subregions
$R_{j_1}^{(\ell+1)}$, $R_{j_2}^{(\ell+1)}$, $R_{j_3}^{(\ell+1)}$ and $R_{j_4}^{(\ell+1)}$,  and then calculating the integral using the Lobatto quadrature rule over these four subregion regions.
The recursive procedure stops when the relative difference of the two approximations at level $\ell$ and $\ell+1$ is smaller than the prescribed
tolerance. We refer the readers to \cite{Berntsen} for more details of the recursive procedure.
Since the region $\Omega_\tau$ is usually thin with small $\tau$ ($\tau=0.1$ is chosen in variety of numerical experiments demonstrated in Section \ref{sec:num_examples}),
the cardinal number of the set $\{ (m,n) | E_{mn}\cap\Omega_\tau\neq\emptyset \} \ll M^2$. In addition,
the recursive adaptive Lobatto quadrature rule for the integral inside $\Omega_\tau$ convergences fast, thus
the adaptive integration only accounts for a small percentage of the overall cost in assembling the stiff matrices.

\subsection{The two-level a posteriori error estimator and mesh refinement}\label{sec:posteriori_err}
Let $\Gamma_h=\{ E_1, E_2, \cdots, E_M \}$  be a mesh over $\Gamma$ with the mesh size $h$, and $\hat\Gamma_h=\{ \hat E_1, \hat E_2, \cdots, \hat E_{2M} \}$
be a uniform refinement of $\Gamma_h$ with the mesh size $h/2$. The corresponding Galerkin solution in the finite-dimensional space $\mathcal{H}_h=\mathcal{U}_h\times\mathcal{V}_h$ 
and  $\hat{\mathcal{H}}_h=\hat{\mathcal{U}}_h\times\hat{\mathcal{V}}_h$
is denoted as $\boldsymbol{\phi}_h$ and $\hat{\boldsymbol{\phi}}_h$, respectively. 
We use a two-level a posteriori error estimator where the exact solution in the numerical error is replaced by the solution $\hat{\boldsymbol{\phi}}_h$ obtained over the uniformly refined mesh $\hat\Gamma_h$.
To this end, we define the first $h-h/2$ based estimators as follows:
\begin{equation}\label{eq:eta}
\eta^{(1)} =   \|\hat\phi^{(1)}_h -\phi^{(1)}_h \|_{H^{1/2}(\Gamma)} \quad  \mbox{and} \quad \eta^{(2)} =  \|\hat\phi^{(2)}_h -\phi^{(2)}_h\|_{H^{-1/2}(\Gamma)},
\end{equation}
in which $\boldsymbol{\phi}_h=[\phi^{(1)}_h, \phi^{(2)}_h]$ and $\hat{\boldsymbol{\phi}}_h=[\hat\phi^{(1)}_h, \hat\phi^{(2)}_h]$.

The estimators $\eta^{(1)}$ and $\eta^{(2)}$ defined above require the computation of the solutions $\boldsymbol{\phi}_h$ and $\hat{\boldsymbol{\phi}}_h$ at both discretization levels,
which could be computationally burdensome. Furthermore, $\hat{\boldsymbol{\phi}}_h$ is expected to be more accurate than $\boldsymbol{\phi}_h$,
thus the latter becomes a temporary result that is not useful once  $\eta^{(1)}$ and $\eta^{(2)}$ are calculated.
In order to reduce the computational cost and avoid such redundancy,
after $\hat{\boldsymbol{\phi}}_h$ is computed, we use $\hat{\boldsymbol{\phi}}_h$ to approximate the solution $\boldsymbol{\phi}_h$
by projecting the refined solution  $\hat{\boldsymbol{\phi}}_h$ over the finite-dimensional space $\mathcal{H}_h$ \cite{Feischl}.
In addition, we localize the estimators by using $h^{1/2}$-weighted $H_1$ and $L_2$ seminorm for $\phi^{(1)}_h$ and $\phi^{(2)}_h$, respectively.
The localization yields error indicator over each element that is computable and can be used to design the mesh refinement strategy.
More precisely, we define the second error estimators $\tilde\eta^{(1)}$ and $\tilde\eta^{(2)}$ as follows:
\begin{equation*}
\tilde\eta^{(1)} =  \left(\sum_{m=1}^{M}  \rho^{(1)}(E_m)\right)^{1/2}\quad  \mbox{and} \quad   \tilde\eta^{(2)} =  \left(\sum_{m=1}^{M}  \rho^{(2)}(E_m)\right)^{1/2},
\end{equation*}
where the error indicators over each element $E_m$ are given by
\begin{equation}\label{eq:rho}
\rho_m^{(1)}(E_m) = |E_m| \cdot \norm {\frac{d}{ds}\left(\hat\phi^{(1)}_h-\Pi_h^{(1)} \hat\phi^{(1)}_h\right) }_{L^2(E_m)}^2 \quad  \mbox{and} \quad
\rho_m^{(2)}(E_m) = |E_m|  \cdot \norm {\hat\phi^{(2)}_h -\Pi_h^{(2)} \hat\phi^{(2)}_h }_{L^2(E_m)}^2,
\end{equation}
and $\Pi_h^{(1)}$ $\left(\Pi_h^{(2)}\right)$ denotes the $L^2$-projection operator from $\hat{\mathcal{U}}_h$ to $\mathcal{U}_h$ ($\hat{\mathcal{V}}_h$ to $\mathcal{V}_h$).
The total error estimator $\tilde\eta$ for $\boldsymbol{\phi}_h$ is defined as
\begin{equation}\label{eq:eta_tilde}
\tilde\eta =  \left( \sum_{m=1}^{M}  \rho^{(1)}(E_m) + \rho^{(2)}(E_m) \right)^{1/2}.
\end{equation}
Recall that $\hat\phi^{(1)}_h$ and $\hat\phi^{(2)}_h$ are the Galerkin approximations of the solution for the transmission problem 
and its normal derivative, we see that $\tilde\eta$ provides a $h^{1/2}$-weighted $H_1$-seminorm estimator for the solution $u$ over the interface $\Gamma$. 
We point out that $\tilde\eta^{(j)}$ ($j=1,2$) alone has been used as a posteriori error estimator for solving elliptic boundary value problems
with Dirichlet or Neumann boundary conditions, and it has been shown that the estimator is efficient and reliable in the sense that the true error for the numerical solution
is bounded below and above by the estimator $\tilde\eta^{(j)}$ \cite{Feischl}.  Here we combine the two together in \eqref{eq:eta_tilde} as the error estimator  for the transmission problem
in the hyperbolic media.  A variety of numerical examples in Section \ref{sec:num_examples}
demonstrate the effectiveness of the proposed error estimator. The theoretical investigation on its robustness remains to be investigated in the future.

The D\"{o}rfler strategy is employed to mark and refine the mesh $\Gamma_h$ using the error estimator $\tilde\eta$.
For given $0<\gamma<1$, we find the minimal set $\mathcal{E}\subset\Gamma_h$ such that
\begin{equation}\label{eq:dorfler}
\gamma \, \tilde\eta^2 \le  \sum_{E\in \mathcal{E} } (\rho^{(1)}(E) + \rho^{(2)}(E)),
\end{equation}
and each marked element in $\mathcal{E}$ is then divided into two sub-elements with equal size. 
The complete adaptive strategy starts from an initial discretization of the interface $\Gamma$.
The calculation of the estimator \eqref{eq:eta_tilde}  and the refinement procedure \eqref{eq:dorfler} is repeated until $\tilde\eta<\sigma$ for certain prescribed tolerance $\sigma$.
This is summarized in the following algorithm:
\begin{algorithm*}
	\caption*{\textbf{Algorithm: The adaptive boundary element method}}
	\begin{algorithmic}
        \State Given tolerance $\sigma$ and the initial mesh size $h$, generate the initial $\Gamma_h^{(0)}$;
		\For {$\ell=0,2,\ldots, L$}
		                \State Let $\hat\Gamma_h^{(\ell)}$ be a uniform refinement of $\Gamma_h^{(\ell)}$;
		                 \State Apply the adaptive numerical integration over $\hat\Gamma_h^{(\ell)}$ to assemble the matrix in \eqref{eq:linear_system};
				\State Compute $\hat{\boldsymbol{\phi}}_h^{(\ell)}$ via  \eqref{eq:linear_system} and compute the estimator $\tilde\eta$ by \eqref{eq:rho} and \eqref{eq:eta_tilde};
			         \If {$\tilde\eta\ge\sigma$}
			                \State  Choose the minimal subset $\mathcal{E}\subset\Gamma_h^{(\ell)}$ such that \eqref{eq:dorfler} is satisfied;
			                \State  Divide each element in $\mathcal{E}$ into two sub-elements to obtain $\Gamma_h^{(\ell+1)}$;
			         \Else
			           \State Stop and return the solution $\hat{\boldsymbol{\phi}}_h^{(\ell)}$;
			         \EndIf
			\EndFor
	\end{algorithmic}
\end{algorithm*}

\section{Numerical examples}\label{sec:num_examples}
We test the accuracy and efficiency of the adaptive boundary element method (BEM) in this section.
Without loss of generality, we consider the point source functions in \eqref{eq:trans_prob} when either $f_1$ or $f_2$ is a Dirac delta function.
Throughout all the examples,  $\tau$ is set as $0.1$ in \eqref{eq:Omega_tau} for the domain $\Omega_\tau$.
The first-order linear element and the zero-order constant element are used to approximate $\phi^{(1)}$ and $\phi^{(2)}$, respectively.
The corresponding numerical solutions returned from the adaptive algorithm are denoted by $\hat\phi_h^{(1)}$ and $\hat\phi_h^{(2)}$,
which are obtained over the refined mesh $\hat\Gamma_h^{(\ell)}$. We also define the relative errors
by letting
$$
 \hat e^{(1)} = \frac{\|\phi^{(1)}-\hat\phi_h^{(1)}\|_{L^2(\Gamma)}}{\|\phi^{(1)}\|_{L^2(\Gamma)}}
 \quad\mbox{and}\quad
 \hat e^{(2)} = \frac{\|\phi^{(2)}-\hat\phi_h^{(2)}\|_{L^2(\Gamma)}}{\|\phi^{(2)}\|_{L^2(\Gamma)}},
$$
where $\phi^{(1)}$ and $\phi^{(2)}$ are reference solutions obtained with  high-order accuracy. 
\\

\noindent\textbf{Example 1} \; Assume that the interface $\Gamma$ is an ellipse with semi-axes $a=2$ and $b=1$ respectively.
The interior domain is a hyperbolic medium while the exterior domain is vacuum.
The permittivity values in $\Omega_1$ and $\Omega_2$ are given by
$\varepsilon_1^{(1)}=1+0.02i$, $\varepsilon_2^{(1)}=-2+0.02i$ and  $\varepsilon_1^{(2)}=\varepsilon_2^{(2)}=1$, respectively. 
The wavenumber $k_0=1$, and the source is located in $\Omega_1$ such that $f_1=-\delta(x)$ and $f_2=0$. 
Due to the smoothness of the interface, we compute the reference solutions $\phi^{(1)}$ and $\phi^{(2)}$ with high-oder accuracy by the Nystr\"{o}m scheme, which is a spectral method
using the trigonometry interpolant over $\Gamma$ \cite{Kress}.
By parameterizing the interface $\Gamma$ with $x_1=a\cos\theta$ and $x_2=b\sin\theta$,
the initial mesh $\Gamma_h^{(0)}$ is generated with grid points $\big\{ (a\cos\theta_m, b\sin\theta_m) \big\}_{m=1}^{M^{(0)}}$, in which $\theta_m = \frac{2(m-1) \pi}{M^{(0)}}$.

\begin{figure}[!htbp]
\begin{center}
\hspace*{-1.2cm}
\includegraphics[height=4cm]{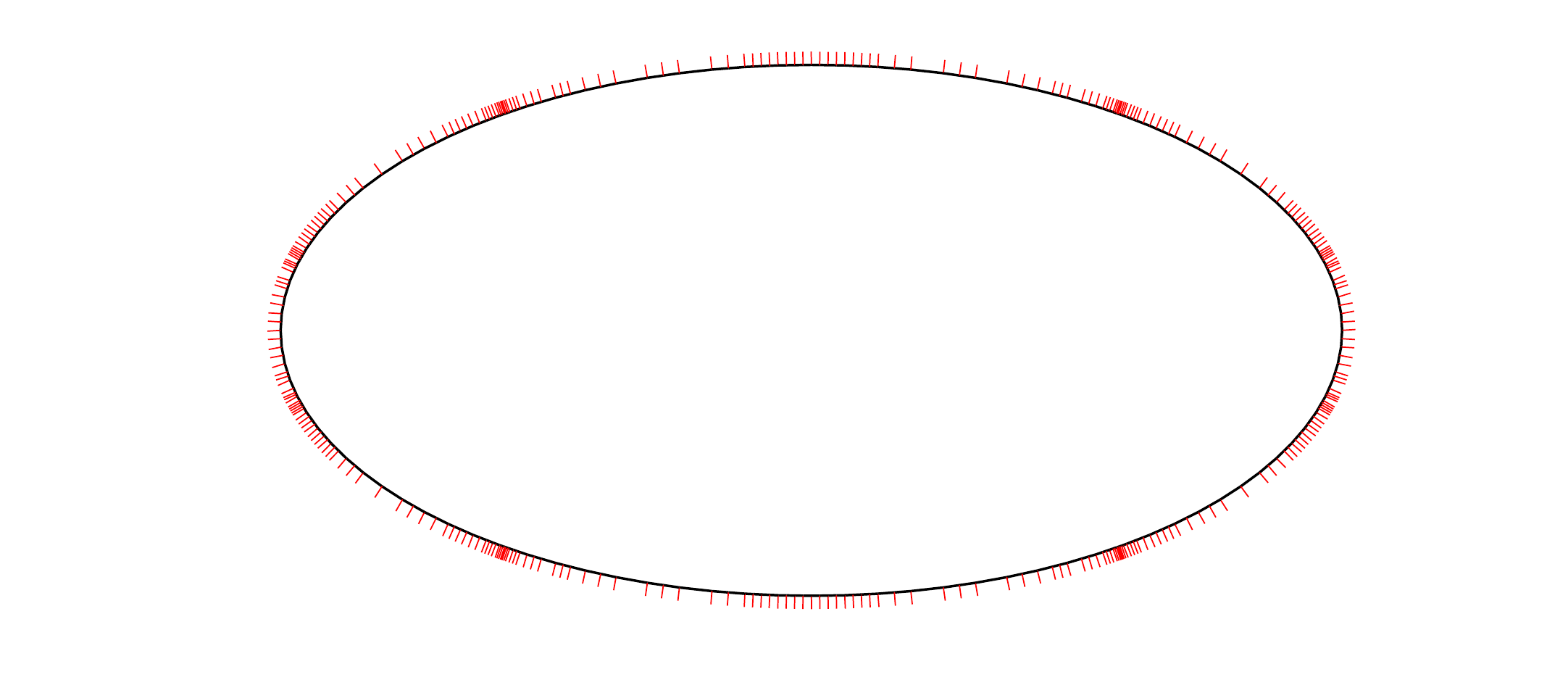}
\hspace*{-1.2cm}
\includegraphics[height=4cm]{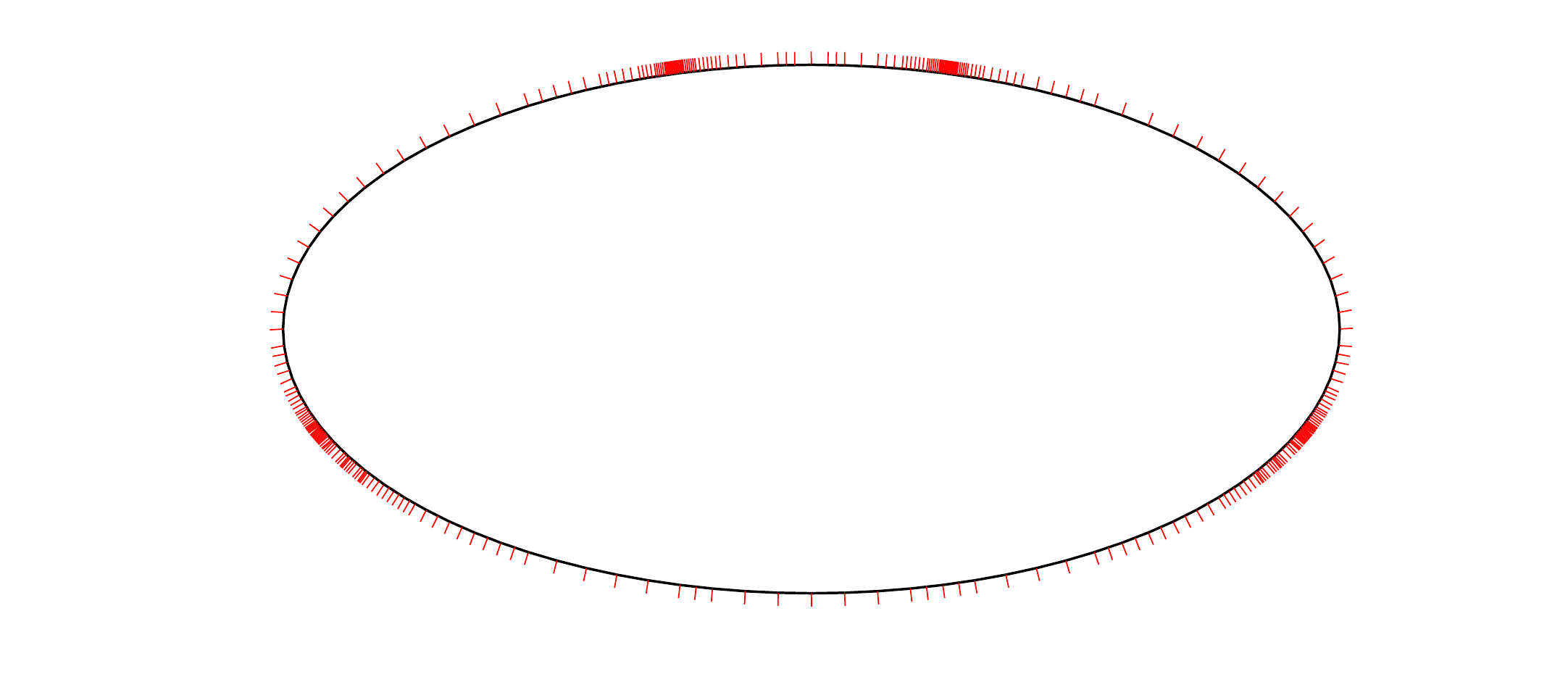}
\caption{The mesh $\Gamma_h^{(4)}$ over the boundary of the ellipse for Example 1 (left) and Example 2 (right).
 }\label{fig:mesh_ex1}
\end{center}
\end{figure}

\begin{figure}[!htbp]
\begin{center}
\includegraphics[height=6cm]{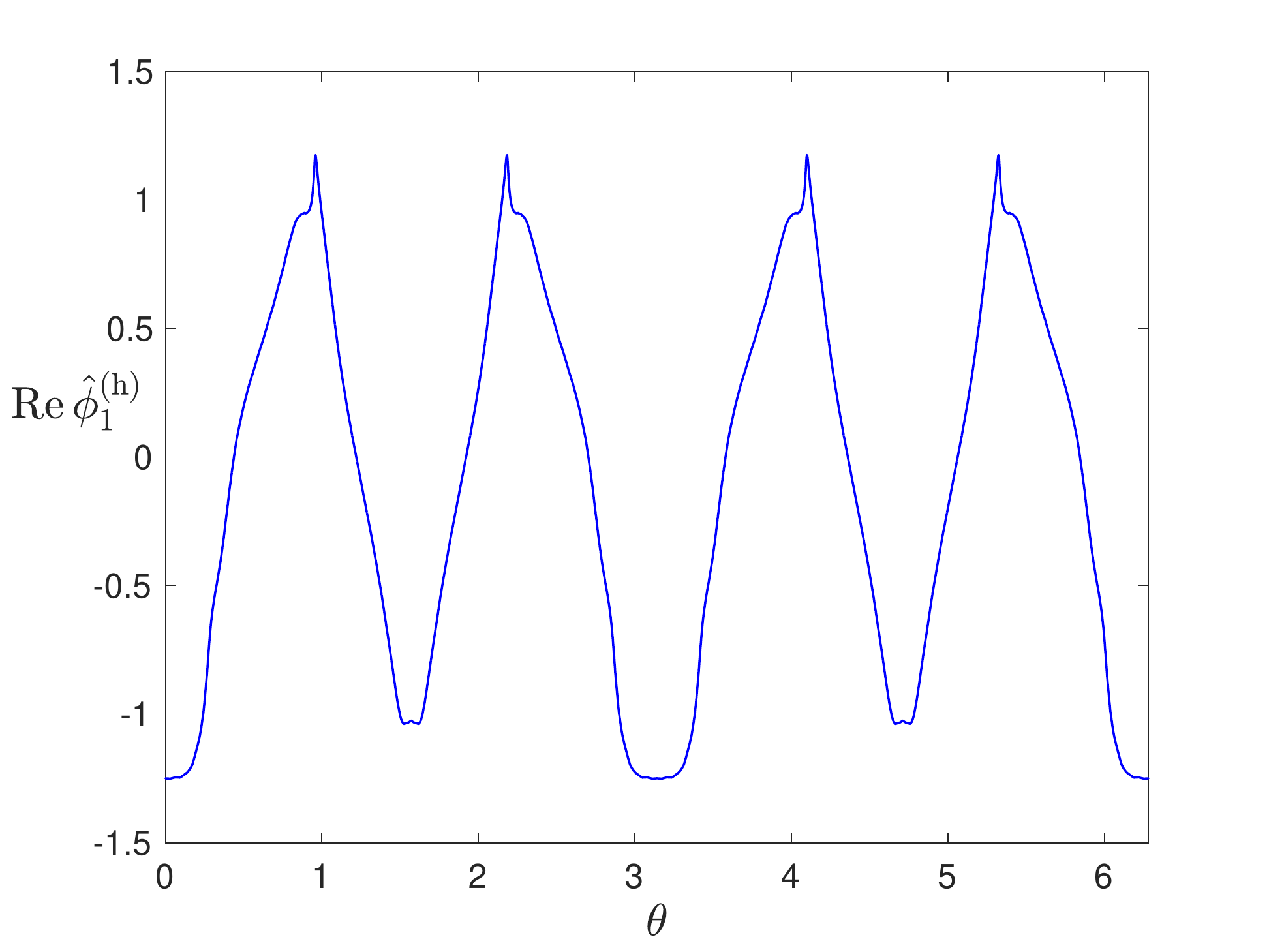}
\includegraphics[height=6cm]{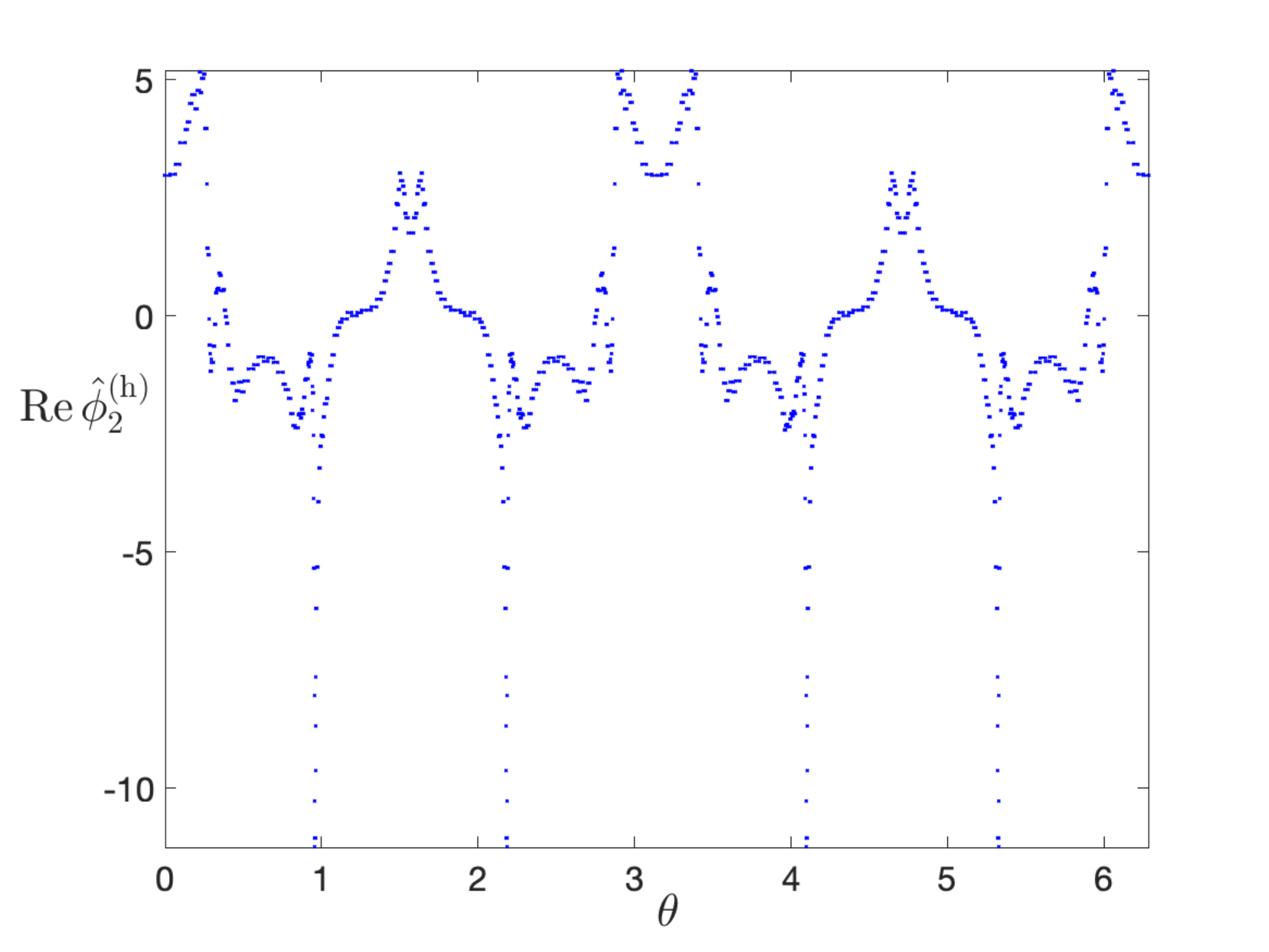}
\caption{The real parts of $\hat\phi_h^{(1)}$ and $\hat\phi_h^{(2)}$ obtained from the adaptive algorithm.
$\theta$ represents the polar angle over the ellipse.
 }\label{fig:sol_ex1}
\end{center}
\end{figure}

\begin{figure}[!htbp]
\begin{center}
\includegraphics[height=6cm]{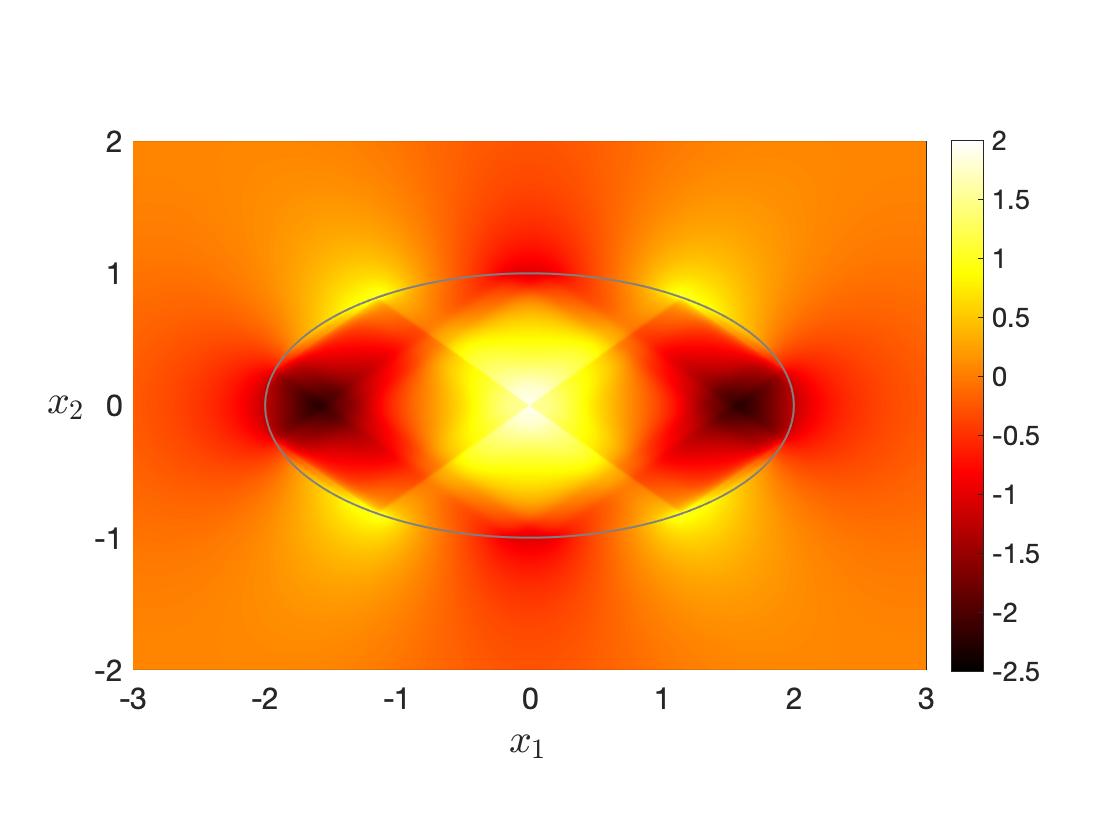}
\includegraphics[height=6cm]{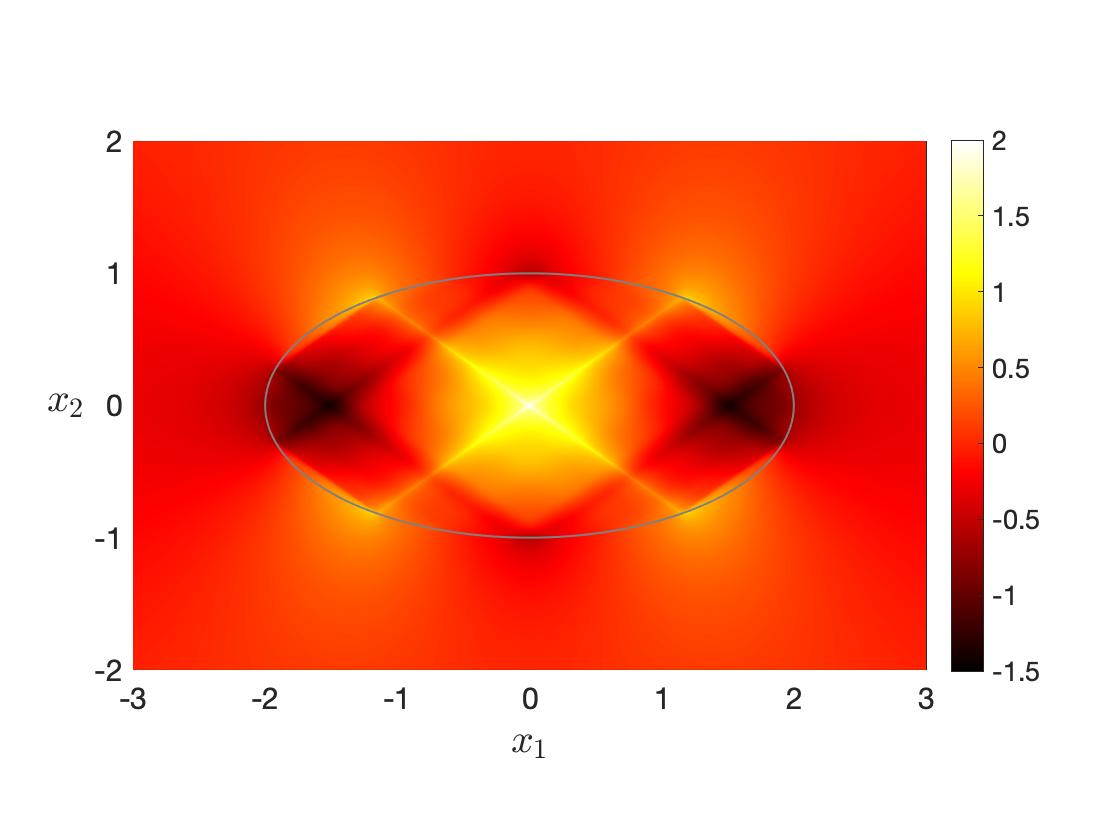}
\caption{Real (left) and imaginary part (right) of the wave field inside and outside the elliptical domain..
 }\label{fig:u_domain_ex1}
\end{center}
\end{figure}

Table \ref{tab:ex1} shows the mesh sizes and the corresponding numerical errors for different levels of refinements when the adaptive procedure is applied,
in which $M^{(\ell)}$ denotes the number of grid points being used in the mesh $\Gamma_h^{(\ell)}$.
 As expected, the numerical error $\hat e^{(1)}$ for $\phi^{(1)}$ is relatively small over the initial mesh  $\Gamma_h^{(0)}$, while
 $\hat e^{(2)}$ is large due to the singular behavior of $\phi^{(2)}$. 
It is observed that the two-level a posteriori error estimator $\tilde\eta$ is effective in identifying the solution singularities and 
a local mesh refinement near the singularities reduces the numerical error significantly after each refinement.
The adaptive procedure terminates with $\ell=4$, and there holds $M^{(4)}=253$ for the mesh $\Gamma_h^{(4)}$;
see Figure \ref{fig:mesh_ex1} (left) for a plot of the mesh $\Gamma_h^{(4)}$.
The numerical solutions $\hat\phi_h^{(1)}$ and $\hat\phi_h^{(2)}$ in the final stage of the adaptive procedure are plotted in Figure \ref{fig:sol_ex1},
which are solved over the mesh $\hat\Gamma_h^{(4)}$.
The corresponding numerical errors are $\hat e^{(1)}=0.0044$ and $\hat e^{(2)}=0.0745$, respectively.
As a comparison, if the mesh is quasi-uniform, one needs a much more refined mesh $\hat\Gamma_h$ with a total number of $2\bar M$ grid points to achieve
a comparable accuracy, in which $\bar M=700$. This is illustrated in the last column of Table \ref{tab:ex1}.
In Figure \ref{fig:u_domain_ex1} we also plot the wave field in the domain $\Omega_1$ and  $\Omega_2$, which are computed via the formulas
\eqref{eq:u_omega1} and \eqref{eq:u_omega2} using the solutions  $\hat\phi_h^{(1)}$ and $\hat\phi_h^{(2)}$ obtained from the adaptive algorithm.
The strong directional propagation of the wave inside the hyperbolic medium and the multiple reflections by the interface $\Gamma$ is clearly seen.

\begin{table}[htbp]
\begin{center}
\caption{The maximum and minimum mesh sizes $h_{\rm{max}}$ and $h_{\rm{min}}$ for $\Gamma_h^{(\ell)}$,
and the relative errors $\hat e^{(1)}$ and $\hat e^{(2)}$ for the adaptive BEM  and non-adaptive BEM. }
\medskip  
\begin{tabular}{||c|ccccc|c||}
\hline
 & \multicolumn{5}{c|}{Adaptive BEM}  & Non-adaptive BEM \\
\hline
 DOF   &  $M^{(0)}=100$   &   $M^{(1)}=125$    &  $M^{(2)}=158$    & $M^{(3)}=202$ & $M^{(4)}=253$  & $\bar M = 700$ \\ 
\hline
\hline    
$h_{\rm{max}}$  & $0.1256$   & $0.1256$  &  $0.1245$ &  $0.1245$  & $0.1234$ & $0.018$ \\
$h_{\rm{min}}$   & $0.0629$   & $0.0349$  &  $0.0174$ &  $0.0087$  & $0.0069$ & $0.009$ \\
\hline
$\hat e^{(1)}$  &  $0.0359$  & $0.0145$ & $0.0057$ &  $0.0047$  & $0.0044$ & $5.28\times10^{-4}$ \\
$\hat e^{(2)}$  &  $0.5685$  &  $0.3804$ & $0.2965$ & $0.1070$  & $0.0745$ & $0.0757$ \\
\hline
\end{tabular}
\label{tab:ex1}
\vspace*{-20pt}
\end{center}
\end{table}

\newpage

\noindent\textbf{Example 2} \; The geometry in this example is the same as in the Example 1, but both media in $\Omega_1$ and $\Omega_2$
are now hyperbolic.  Their permittivity values are given by $\varepsilon_1^{(1)}=-1+0.02i$, $\varepsilon_2^{(1)}=1+0.02i$,  and $\varepsilon_1^{(2)}=-4+0.05i$, $\varepsilon_2^{(2)}=1+0.05i$.
Assume that the source is located in exterior domain such that $f_1=0$ and $f_2=-\delta(x-x_0)$, in which the source location is $x_0=(0,2)$.

\begin{figure}[!htbp]
\begin{center}
\includegraphics[height=6cm]{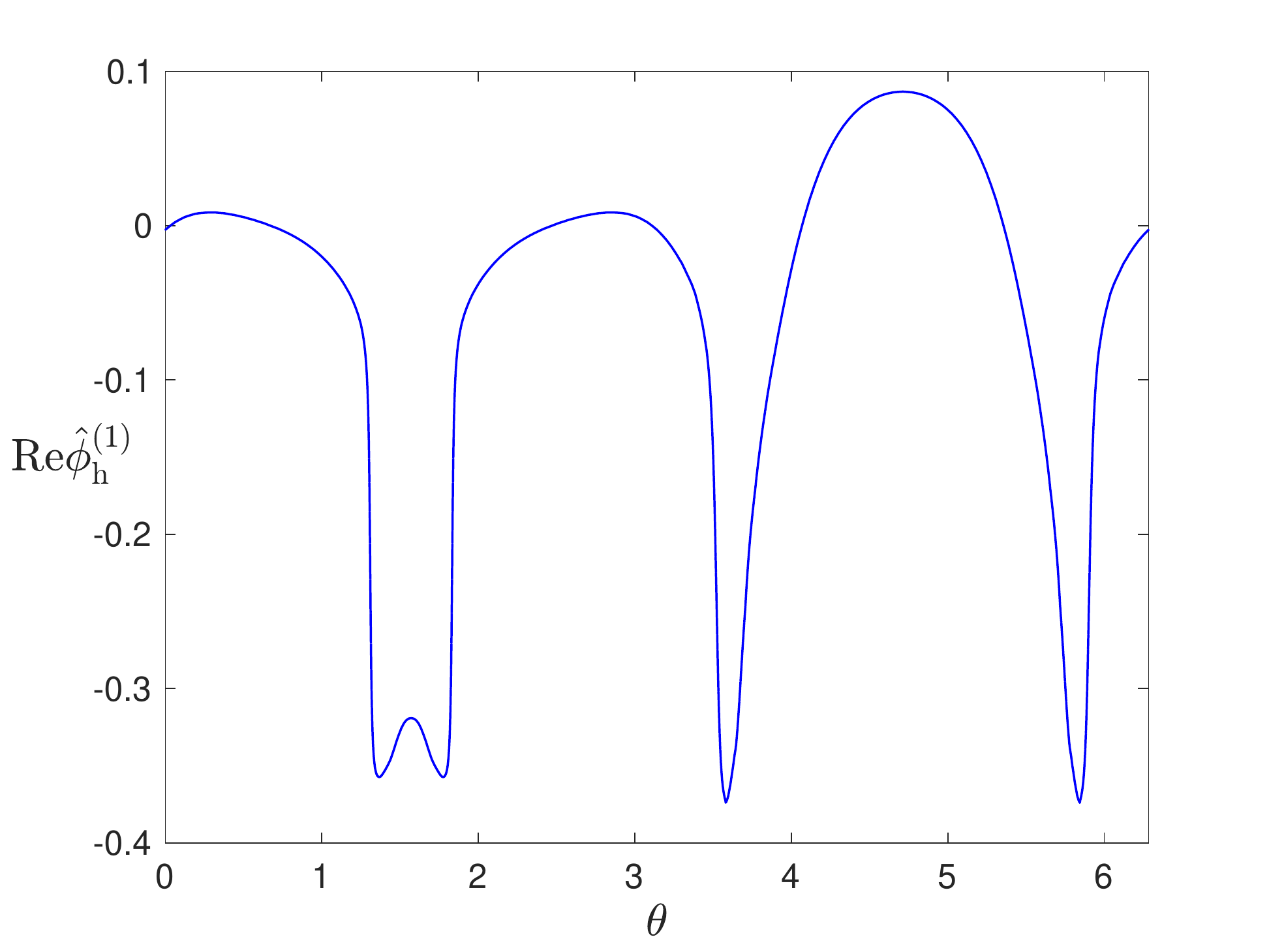}
\includegraphics[height=6cm]{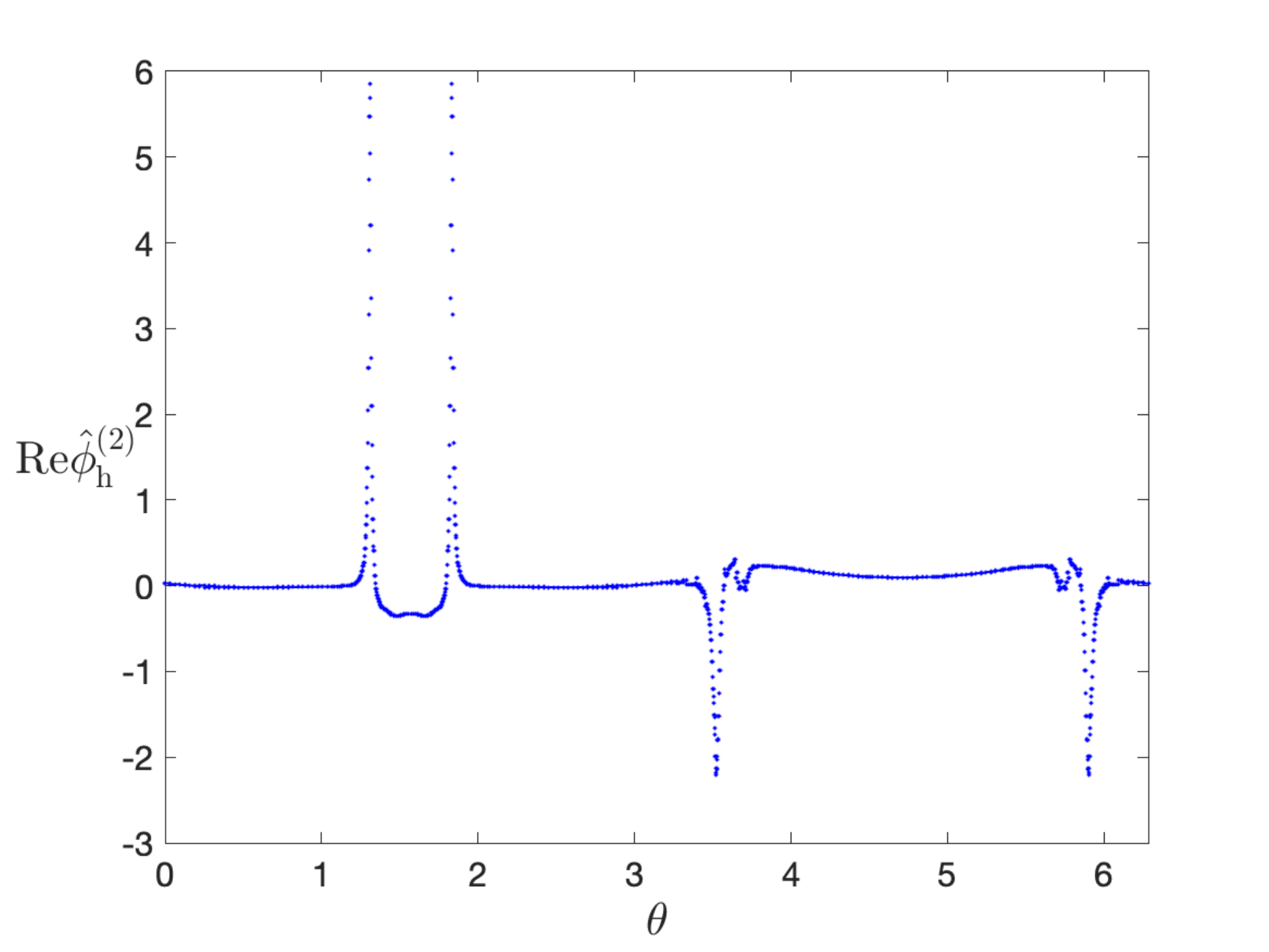}
\caption{The real parts of $\hat\phi_h^{(1)}$ and $\hat\phi_h^{(2)}$ obtained from the adaptive algorithm.
$\theta$ represents the polar angle over the ellipse.
 }\label{fig:sol_ex2}
\end{center}
\end{figure}

\begin{figure}[!htbp]
\begin{center}
\includegraphics[height=6cm]{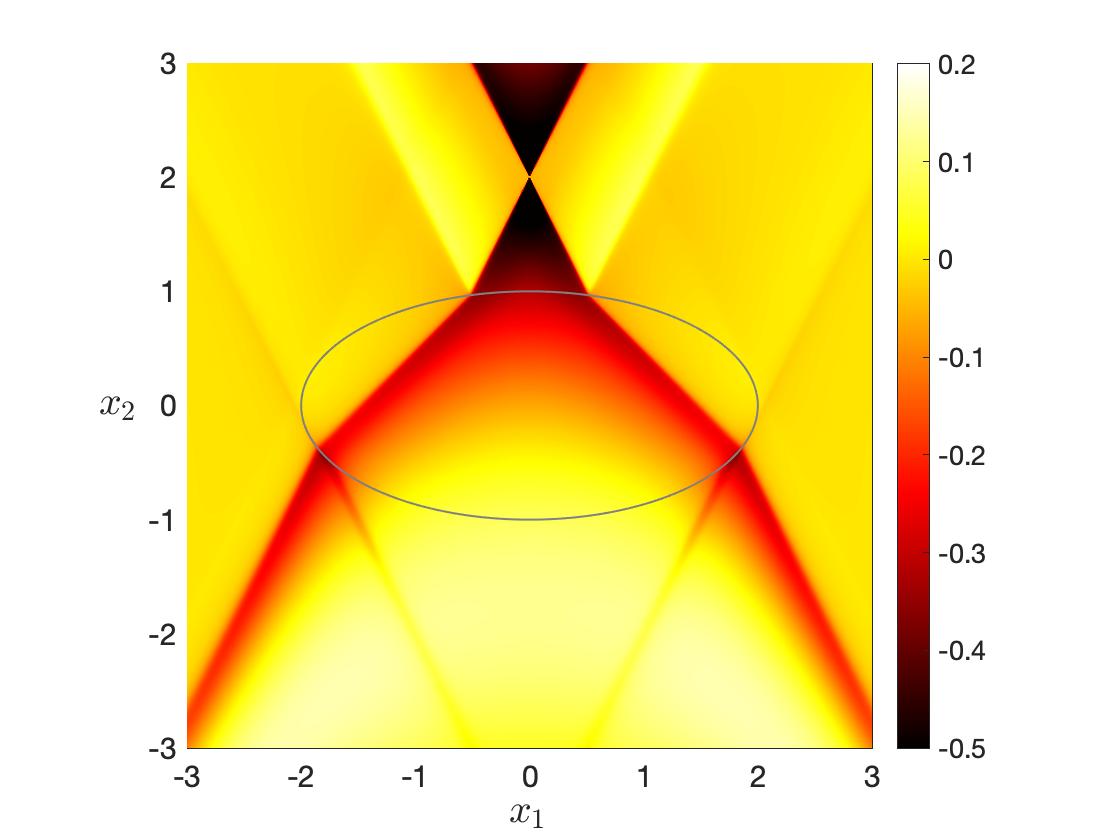}
\includegraphics[height=6cm]{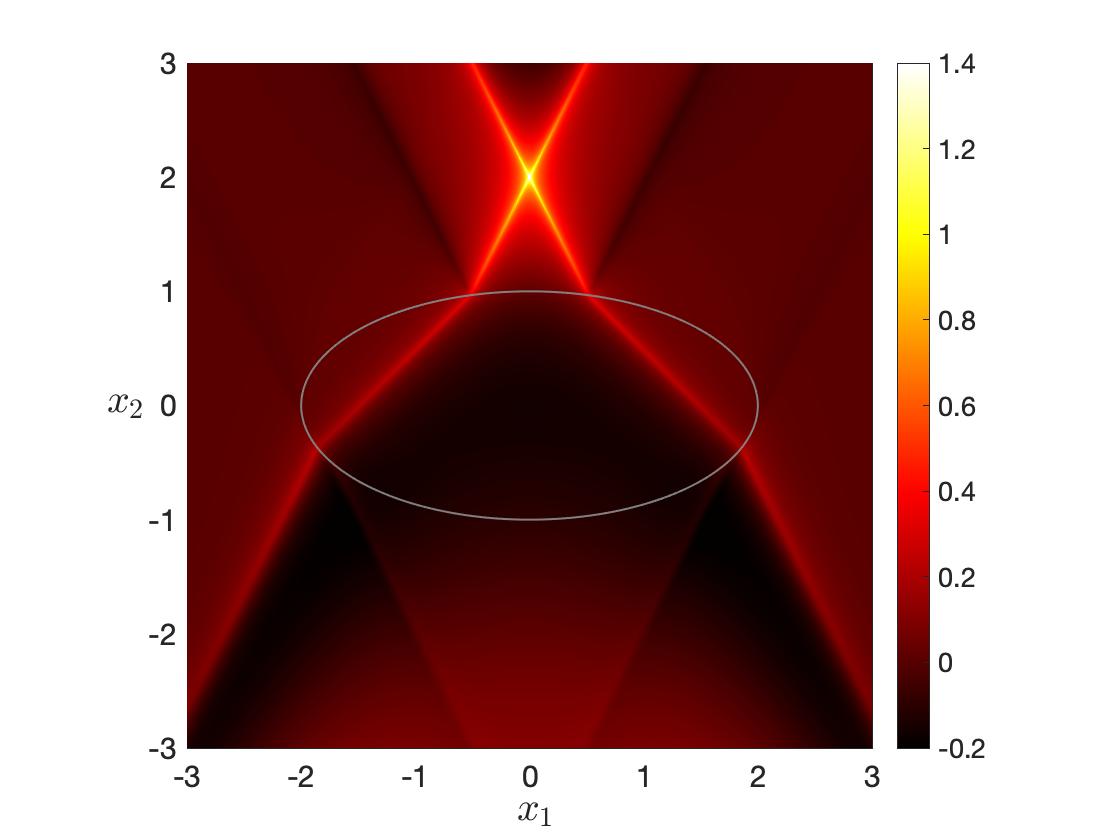}
\caption{Real (left) and imaginary part (right) of the wave field inside and outside the elliptical domain.
 }\label{fig:u_domain_ex2}
\end{center}
\end{figure}

The adaptive procedure also terminates with $\ell=4$ and the mesh $\Gamma_h^{(4)}$ attains a total of $249$ grids points as shown in Figure \ref{fig:mesh_ex1} (right).
The final numerical solutions $\hat\phi_h^{(1)}$ and $\hat\phi_h^{(2)}$ are plotted in Figure \ref{fig:sol_ex2},
and the corresponding numerical errors are $\hat e^{(1)}=0.0034$ and $\hat e^{(2)}=0.0616$, respectively.
For completeness, we collect all the mesh sizes and the corresponding numerical errors for different levels of refinements in Table \ref{tab:ex2}.
The computation with a quasi-uniform mesh to achieve a comparable accuracy is also shown. 
We see that the number of degrees of freedom is reduced by about 4 times when the adaptive procedure is applied.
Figure \ref{fig:u_domain_ex2} demonstrates the wave field in the domain $\Omega_1$ and  $\Omega_2$.
The wave is strongly directional while penetrating through the interior hyperbolic medium and being reflected at the interface. \\

\begin{table}[htbp]
\begin{center}
\caption{The maximum and minimum mesh sizes $h_{\rm{max}}$ and $h_{\rm{min}}$ for $\Gamma_h^{(\ell)}$,
and the relative errors $\hat e^{(1)}$ and $\hat e^{(2)}$ for the adaptive BEM  and non-adaptive BEM. }
\medskip  
\begin{tabular}{||c|ccccc|c||}
\hline
 & \multicolumn{5}{c|}{Adaptive BEM}  & Non-adaptive BEM \\
\hline
 DOF    &  $M^{(0)}=100$   &   $M^{(1)}=125$    &  $M^{(2)}=157$    & $M^{(3)}=199$  & $M^{(4)}=249$ & $\bar M=1000$ \\ 
\hline
\hline 
$h_{\rm{max}}$  & $0.1256$   & $0.1256$  &  $0.1256$ & $0.1256$ & $0.1256$ & $0.0126$ \\
$h_{\rm{min}}$   & $0.0629$   & $0.0364$  &  $0.0182$ & $0.0092$ & $0.0046$ & $0.0063$\\
\hline
$\hat e^{(1)}$  &  $0.0576$  &  $0.0154$ & $0.0042$ &  $0.0034$   &  $0.0034$ & $2.16\times10^{-4}$ \\
$\hat e^{(2)}$  &  $0.6957$  &  $0.4197$ & $0.1881$ &  $0.0928$   &  $0.0616$ & $0.0639$  \\ 
\hline
\end{tabular}
\label{tab:ex2}
\vspace*{-20pt}
\end{center}
\end{table}

\begin{figure}[!htbp]
\begin{center}
\includegraphics[height=3.5cm]{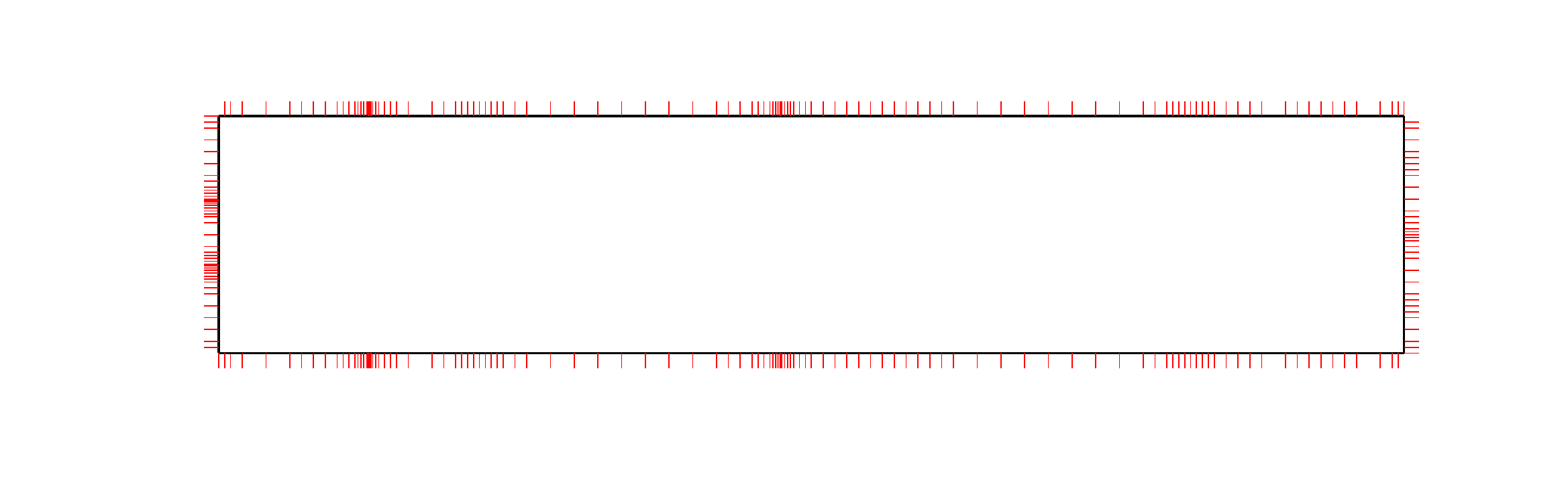}
\includegraphics[height=3.5cm]{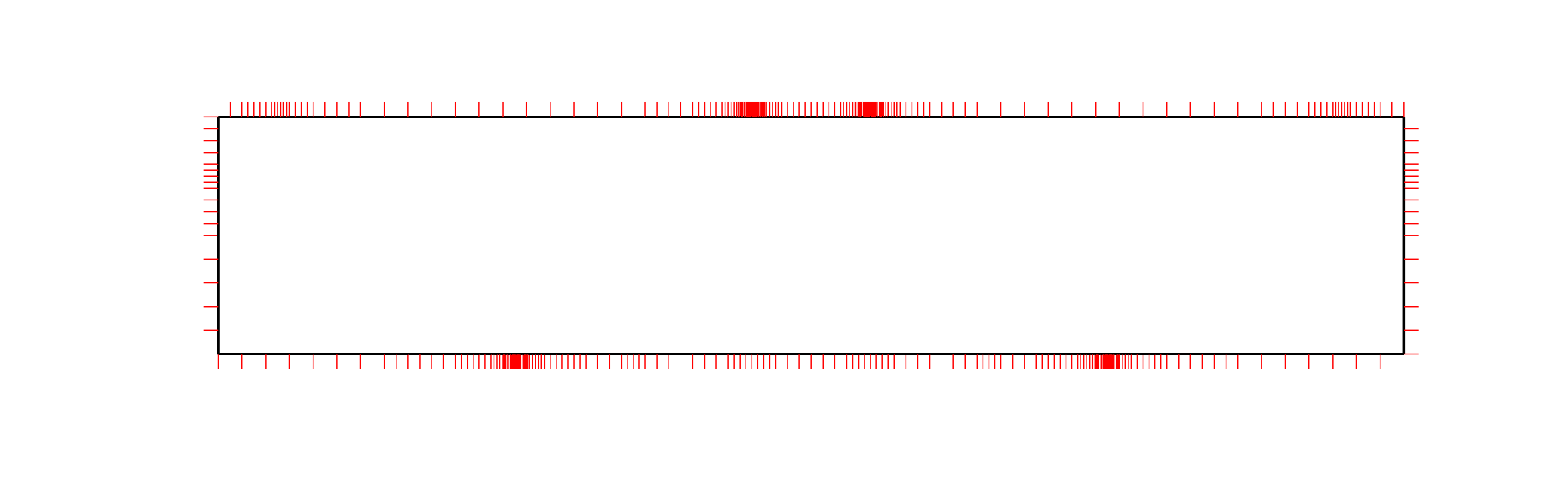}
\caption{The mesh $\Gamma_h^{(3)}$ (top)  and $\Gamma_h^{(4)}$ (bottom) over the boundary of the rectangle for Example 3 and Example 4 respectively.
 }\label{fig:mesh_ex3}
\end{center}
\end{figure}

\begin{figure}[!htbp]
\begin{center}
\includegraphics[height=5cm]{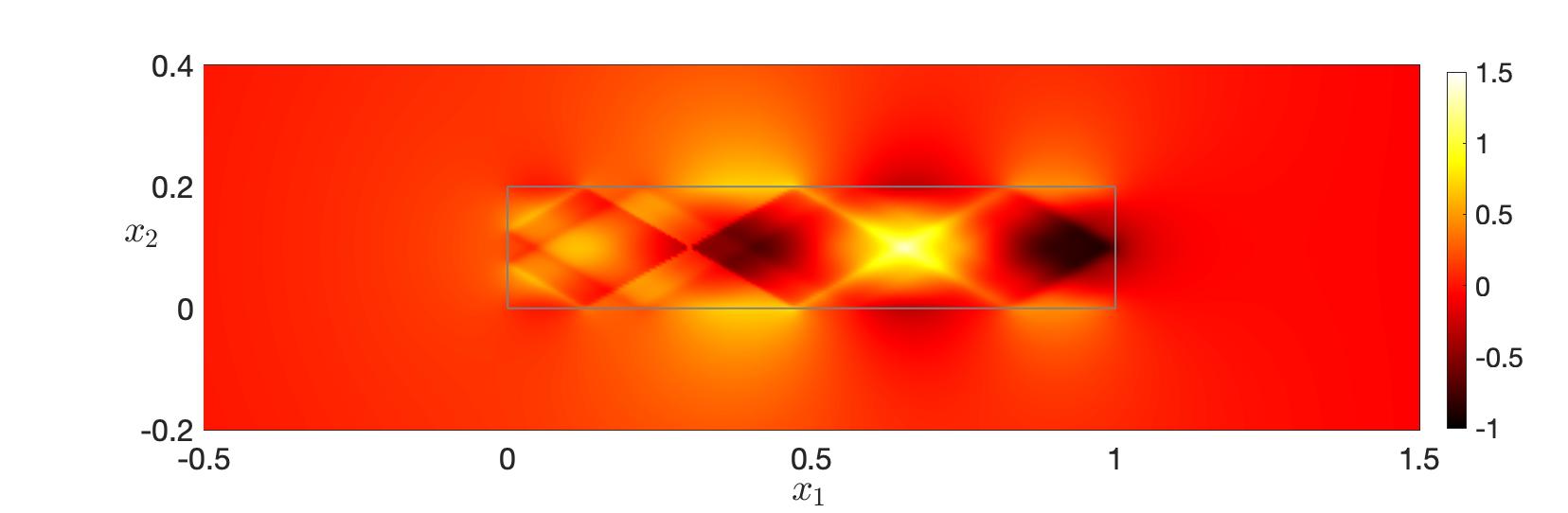}\\
\hspace*{0.015cm}
\includegraphics[height=5cm]{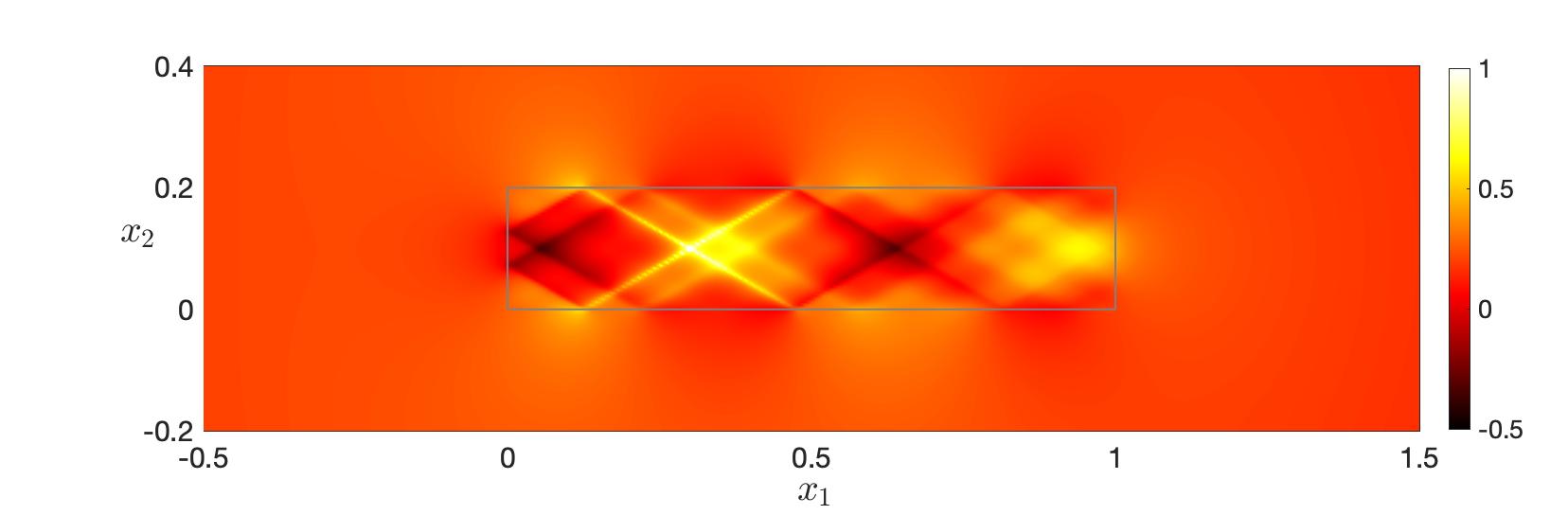}
\caption{Real (top) and imaginary part (bottom) of the wave field inside and outside the rectangle.
 }\label{fig:u_domain_ex3}
\end{center}
\end{figure}

\newpage
\noindent\textbf{Example 3} \; In this example, we consider a rectangular domain $\Omega_1= (0, 1) \times (0, 0.2)$ 
with the permittivity values  $\varepsilon_1^{(1)}=1+0.02i$, $\varepsilon_2^{(1)}=-3+0.1i$. The exterior domain $\Omega_2$ is occupied by vacuum. 
When the source is located at $(0.3, 0.1)$ and the wavenumber $k_0=1$, 
the mesh sizes and the corresponding numerical errors for the adaptive algorithm are shown in Table \ref{tab:ex3}.
Here we obtain the reference solutions $\phi^{(1)}$ and $\phi^{(2)}$ with high-oder accuracy by using a uniform mesh with large degree-of-freedom.
The final mesh $\Gamma_h^{(3)}$ attains a total of $238$ grids points as shown in  Figure \ref{fig:mesh_ex3} (top).
The numerical errors for the solutions obtained from $\hat\Gamma_h^{(3)}$  are $\hat e^{(1)}=0.0011$ and $\hat e^{(2)}=0.0722$, respectively.
In contrast, the same level of accuracy is obtained by a quasi-uniform mesh $\hat\Gamma_h$ with a total number of $2\bar M$ grid points, in which $\bar M=672$.
Figure \ref{fig:u_domain_ex3} plots the wave field in the domain $\Omega_1$ and  $\Omega_2$.
It is observed that multiple reflections by the interface induce strong wave interactions inside $\Omega_1$. \\

\begin{table}[htbp]
\begin{center}
\caption{The maximum and minimum mesh sizes $h_{\rm{max}}$ and $h_{\rm{min}}$ for $\Gamma_h^{(\ell)}$,
and the relative errors $\hat e^{(1)}$ and $\hat e^{(2)}$ for the adaptive BEM  and non-adaptive BEM. }
\medskip  
\begin{tabular}{||c|cccc|c||}
\hline
 & \multicolumn{4}{c|}{Adaptive BEM}  & Non-adaptive BEM \\
\hline
     &  $M^{(0)}=120$   &   $M^{(1)}=150$    &  $M^{(2)}=190$    & $M^{(3)}=238$ & $\bar M=672$\\ 
\hline
\hline
$h_{\rm{max}}$  & $0.02$   & $0.02$  & $0.02$ & $0.02$  & 0.0036 \\
$h_{\rm{min}}$   & $0.02$   & $0.01$  & $0.0050$ & $0.0025$  & 0.0036 \\
\hline
$\hat e^{(1)}$  &  $0.0333$  & $0.0078$ & $0.0023$ &  $0.0011$  & $3.76\times10^{-4}$ \\
$\hat e^{(2)}$  &  $0.5684$  & $0.2860$ & $0.1268$ &  $0.0722$  & 0.0714\\
\hline
\end{tabular}
\label{tab:ex3}
\vspace*{-20pt}
\end{center}
\end{table}

\newpage

\begin{figure}[!htbp]
\begin{center}
\hspace*{-0.5cm}
\includegraphics[height=5cm]{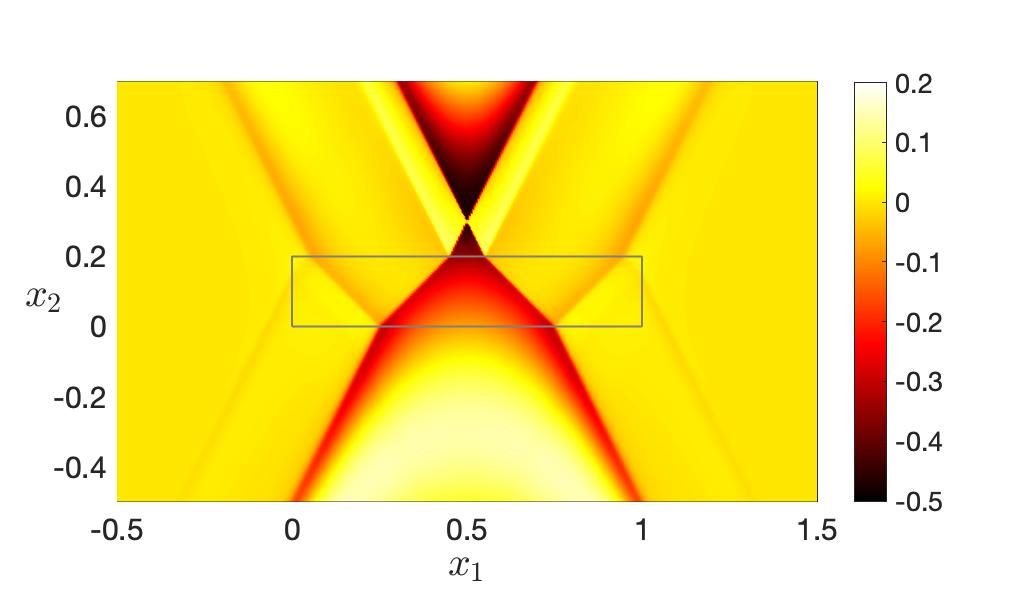}
\hspace*{-0.5cm}
\includegraphics[height=5cm]{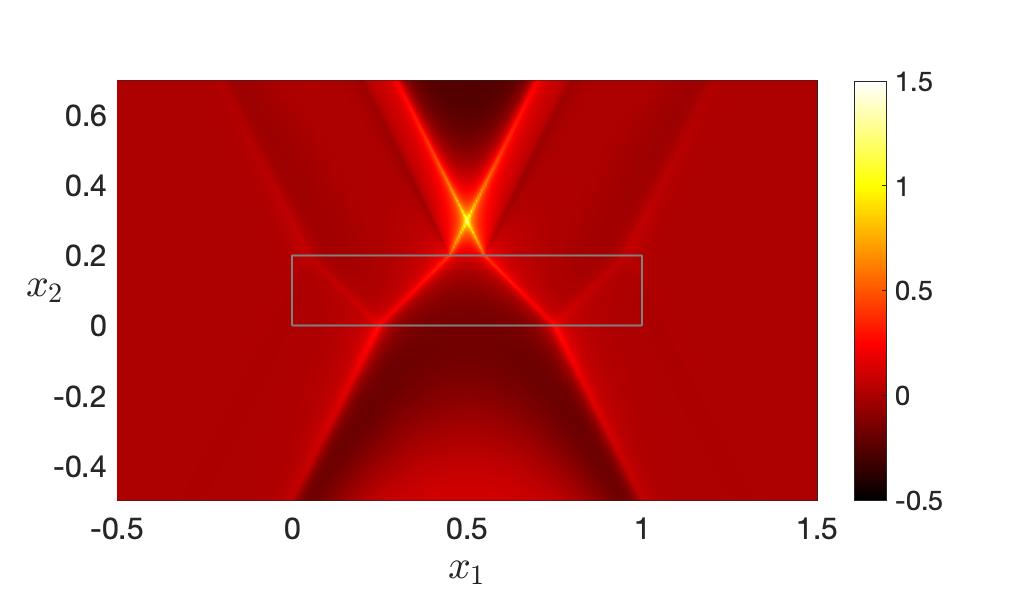}
\caption{Real (left) and imaginary part (right) of the wave field inside and outside the rectangle.
 }\label{fig:u_domain_ex4}
\end{center}
\end{figure}

\noindent\textbf{Example 4} \; We consider the same geometry as in Example 3 while assigning the permittivity values 
as $\varepsilon_1^{(1)}=-1+0.02i$, $\varepsilon_2^{(1)}=1+0.02i$ and $\varepsilon_1^{(2)}=-4+0.05i$, $\varepsilon_2^{(2)}=1+0.05i$
for the interior and exterior domain, respectively. The wavenumber $k_0=2\pi$, and
the source is located in exterior domain with $f_1=0$ and $f_2=-\delta(x-x_0)$, where the source location is $x_0=(0.5, 0.3)$.
The adaptive boundary element still successfully leads to a reduction of numerical errors to the desired accuracy when the mesh refinement is performed locally,
as demonstrated by Table \ref{tab:ex4}. Due to the singularity of the solution, one needs a quasi-uniform mesh with the mesh size $h_{max}=0.0013$ to achieve
an accuracy that is obtained by an adaptive mesh with $h_{max}=0.02$. The number of degrees of freedom for the former is more than 5 times higher than the latter.
The final mesh in the adaptive procedure and the corresponding wave field in the domain are shown in Figures \ref{fig:mesh_ex3} (right) and \ref{fig:u_domain_ex4},
respectively.

\begin{table}[htbp]
\begin{center}
\caption{The maximum and minimum mesh sizes $h_{\rm{max}}$ and $h_{\rm{min}}$ for $\Gamma_h^{(\ell)}$,
and the relative errors $\hat e^{(1)}$ and $\hat e^{(2)}$ for the adaptive BEM  and non-adaptive BEM. }
\medskip  
\begin{tabular}{||c|ccccc|c||}
\hline
 & \multicolumn{5}{c|}{Adaptive BEM}  & Non-adaptive BEM \\
\hline
     &  $M^{(0)}=120$   &   $M^{(1)}=150$    &  $M^{(2)}=188$    & $M^{(3)}=235$  & $M^{(4)}=294$  & $\bar M=1920$ \\ 
\hline
\hline
$h_{\rm{max}}$  & $0.02$   & $0.02$  & $0.02$ & $0.02$  & $0.02$ & $0.0013$ \\
$h_{\rm{min}}$   & $0.02$   & $0.01$  & $0.005$ & $0.0025$  & $1.25 \times 10^{-3}$  & $0.0012$\\
\hline
$\hat e^{(1)}$  &  $0.1044$  & $0.0365$ & $0.0093$ &  $0.0015$  & $8.31 \times 10^{-4}$ & $7.79 \times 10^{-4}$\\
$\hat e^{(2)}$  &  $0.7599$  & $0.6673$ & $0.4123$ &  $0.1785$  & $0.0812$  & $0.0804$   \\
\hline
\end{tabular}
\label{tab:ex4}
\vspace*{-20pt}
\end{center}
\end{table}

\newpage

\begin{figure}[!htbp]
\begin{center}
\includegraphics[height=5cm]{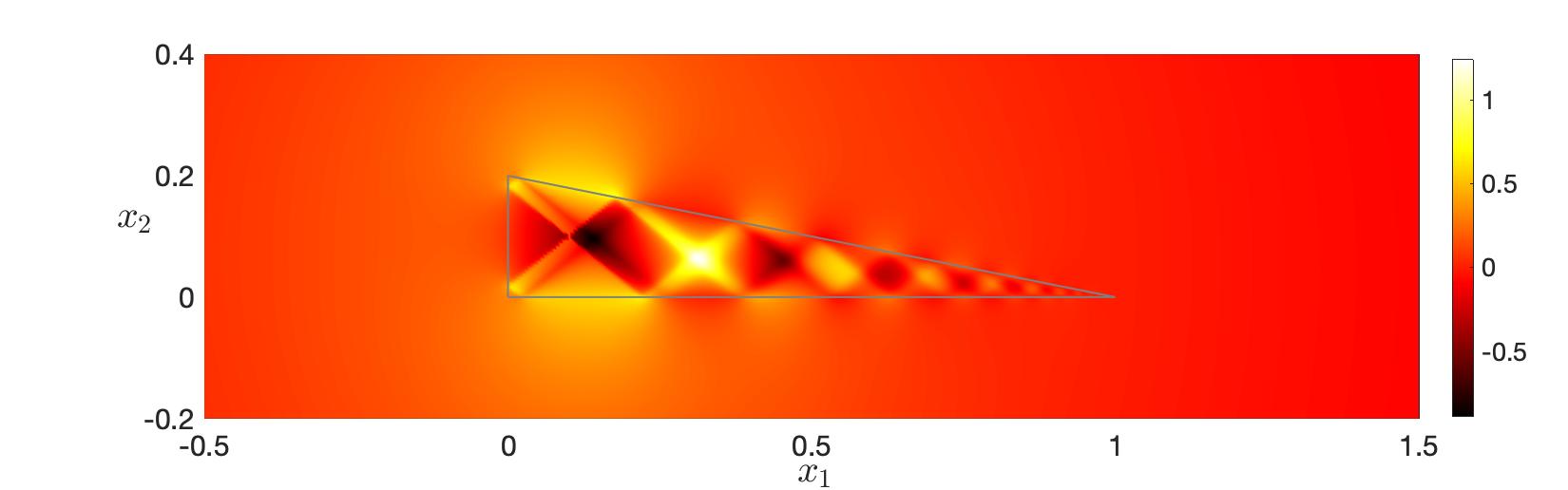}\\
\hspace*{0.015cm}
\includegraphics[height=5cm]{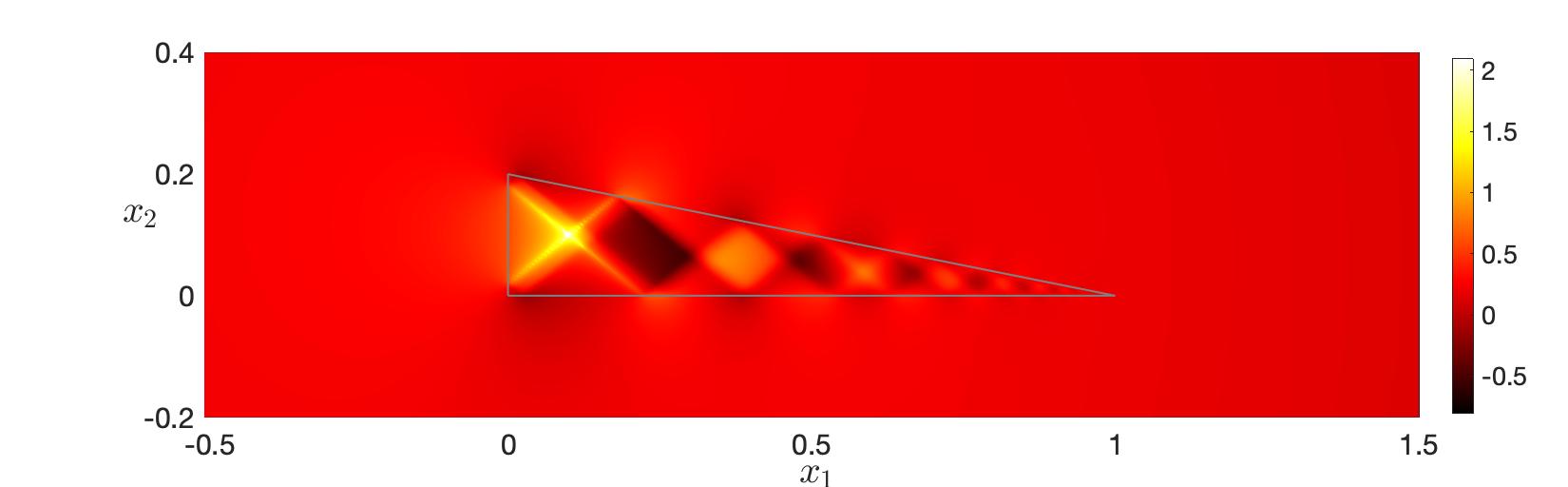}
\vspace*{-5pt}
\caption{Real (top) and imaginary part (bottom) of the wave field inside and outside the wedge.
 }\label{fig:u_domain_ex5}
\end{center}
\end{figure}

\noindent\textbf{Example 5}  \; In the final example, the interior domain $\Omega_1$ attains a wedge shape and it is placed in the vacuum. 
Such geometry has been be used for super-focusing of the electromagnetic waves near the sharp tip \cite{Nikitin}. 
The permittivity values for the hyperbolic medium are
$\varepsilon_1^{(1)}=2$, $\varepsilon_2^{(1)}=-3+0.03i$.  The source is located in $\Omega_1$ so that 
$f_1=-\delta(x-x_0)$ with $x_0=(0.1, 0.1)$ and $f_2=0$, and the wavenumber $k_0=1$.
The adaptive algorithm terminates after three mesh refinements, 
which yield the numerical solution with an accuracy of $\hat e^{(1)}=0.0011$ and $\hat e^{(2)}=0.0722$ (see Table \ref{tab:ex5}).
As demonstrated in Figure \ref{fig:u_domain_ex5}, the wave propagates toward the tip while being reflected by the boundary of the wedge.

\begin{table}[htbp]
\begin{center}
\caption{The maximum and minimum mesh sizes $h_{\rm{max}}$ and $h_{\rm{min}}$ for $\Gamma_h^{(\ell)}$,
and the relative errors  $\hat e^{(1)}$ and $\hat e^{(2)}$ for the adaptive BEM  and non-adaptive BEM. }
\medskip  
\begin{tabular}{||c|cccc|c||}
\hline
 & \multicolumn{4}{c|}{Adaptive BEM}  & Non-adaptive BEM \\
\hline
     &  $M^{(0)}=220$   &   $M^{(1)}=275$    &  $M^{(2)}=344$    & $M^{(3)}=431$  & $\bar M=1210$  \\ 
\hline
\hline
$h_{\rm{max}}$  & $0.0102$   & $0.0102$  & $0.0102$ & $0.0102$  & $0.0019$ \\
$h_{\rm{min}}$   & $0.01$   & $0.005$  & $0.0025$  & $1.25 \times 10^{-3}$ & $0.0018$\\
\hline
$\hat e^{(1)}$  &  $0.0219$  & $0.0110$ & $0.0051$ & $0.0030$ & $0.0062$ \\
$\hat e^{(2)}$  &  $0.4025$  & $0.2515$ & $0.1126$ & $0.0658$  & $0.0696$ \\
\hline
\end{tabular}
\label{tab:ex5}
\vspace*{-20pt}
\end{center}
\end{table}

\begin{figure}[!htbp]
\begin{center}
\includegraphics[height=3cm]{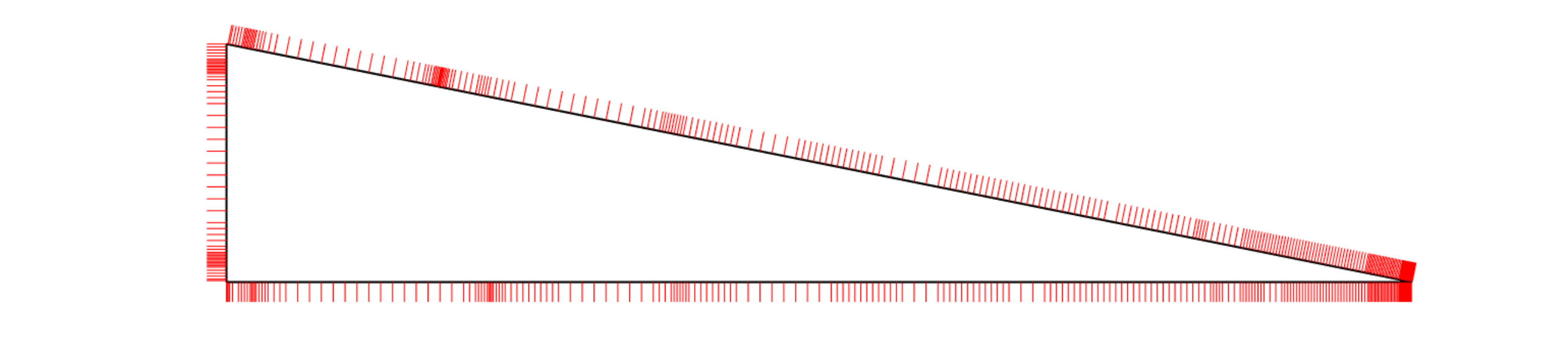}

\caption{The mesh $\Gamma_h^{(3)}$ over the boundary of the wedge.
 }\label{fig:mesh_ex5}
\end{center}
\end{figure}

\section{Discussions}
An adaptive  boundary element method was presented in this paper for solving the transmission problem with hyperbolic metamaterials.
Compared to the discretization with a quasi-uniform mesh, the adaptive approach is able to resolve the singular behavior of the solution with
local mesh refinement, which reduces the number of the degrees of freedom and the overall computational cost significantly.
There are several theoretical and computational issues to be explored along this direction.
Although numerical examples show the efficacy of the error estimator and the accuracy of the adaptive procedure, the robust analysis of the estimator
and the convergence analysis of the algorithm have not been carried out yet. As pointed out previously, this work is mainly on the demonstration of the adaptive algorithm for the two dimensional problems,
its application in three dimensions has not been explored. It is expected that the adaptive algorithm would offer even larger reduction on the computational cost for 3D simulations. 
However, the computation becomes more challenging for mesh refinement and numerical integration in 3D.
Note that for the full Maxwell's equations, the Dyadic Green's functions used in the integral formulation is highly anisotropic with coexistence of the cone-like pattern 
due to emission of the extraordinary TM-polarized waves and elliptical pattern due to emission of ordinary TE-polarized waves \cite{Potemkin}.
Finally, efficient integral equation methods for hyperbolic metamaterials with unbounded domains (e.g. layered media) and other settings of practical interest remain to investigated.
The propagation nature of waves with arbitrarily large wave vectors inside the propagating cone requires new treatement in developing the computational algorithm.

\bibliography{references}

\end{document}